\documentclass[journal]{IEEEtran}
\ifCLASSINFOpdf
\else
\fi
\usepackage[comma,numbers,square,sort&compress]{natbib}
\usepackage{algorithm} %format of the algorithm
\usepackage{algorithmic} %format of the algorithm

\usepackage{multirow} %multirow for format of table
\usepackage{xcolor}
\usepackage{setspace}
 \usepackage{graphicx}
\usepackage{amsmath}
\newtheorem{theorem}{\textbf{Theorem}}
\newtheorem{lemma}{\textbf{Lemma}}
\newtheorem{assumption}{\textbf{Assumption}}
\newtheorem{corollary}{\textbf{Corollary}}
\newtheorem{remark}{\textbf{Remark}}
\newtheorem{definition}{\textbf{Definition}}

\usepackage{multirow}
\usepackage{url}
\usepackage{subfigure}
\usepackage{caption}
\usepackage{geometry}
\usepackage{fancyhdr}
\usepackage{amsmath}
\usepackage{multirow}
\usepackage{amssymb }
\usepackage{color}
\usepackage{graphics} % for pdf,  bitmapped graphics files
\usepackage{graphicx}
\begin{document}
%\begin{spacing}{1.5}
\title{\LARGE \bf
%On Link Selection and Controllability Robustness for Networks: Complexity Analysis
%On Link Selection and Robustness for Network Controllability: Complexity Analysis
%Observability Robustness under Sensor Failures: \\~ Complexities and algorithms
Observability Robustness under Sensor Failures: { \\~ a Computational Perspective}
%On Observability Robustness under Sensor Failures
%Observability Robustness under Sensor Failures: \\~ Complexities Analysis
}
% %\author{Yuan Zhang % <-this % stops a space
%%\thanks{*This work was supported in part by the NNSFC under Grant 61573209 and 61733008. }% <-this % stops a space
\author{Yuan Zhang, Yuanqing Xia, and Kun Liu
\thanks{The authors are with the School of Automation, Beijing Institute of Technology, Beijing, China. {Partial work of this paper was done when the first author was in Tsinghua University (see the preliminary version of this paper \cite{zhang2018arxiv})}. %(emails: zhangyuan14@bit.edu.cn, xia_yuanqing@bit.edu.cn,	kunliubit@bit.edu.cn)
       {(email: {\tt\small \{zhangyuan14,xia\_yuanqing,kunliubit\}@bit.edu.cn}).}}}
\maketitle

% in the abstract or keywords. %malicious
\begin{abstract} %The problem of determining the minimal number of sensors whose removal destroys observability of a linear time invariant system is studied. This problem is closely related to the ability of unique state reconstruction of a system under adversarial sensor attacks, and the dual of it is the inverse to the recently studied minimal controllability problems. It is proven that this problem is NP-hard both for a numerically specific system, and for a generic system whose nonzero entries of its system matrices are unknown but can take values freely (also called structured system). Two polynomial time algorithms are provided to solve this problem, respectively, on a numerical system with bounded maximum geometric multiplicities, and on a structured system with bounded matching deficiencies, which are often met by practical engineering systems.  The proposed algorithms can be easily extended to the case where each sensor has a non-negative cost. Numerical experiments show that the structured system based algorithm could be alternative when the exact values of system matrices are not accessible.

This paper studies the robustness of observability of a linear time-invariant system under sensor failures {from a computational perspective}. To be precise, the problem of determining the minimum number of sensors whose removal can destroy system observability, {as well as the problem of determining the minimum number of state variables that need to be prevented from being directly measured by the existing sensors to destroy observability,} is investigated. The first one is closely related to the ability of unique state reconstruction of a system under adversarial sensor attacks, and the dual of both problems are in the {opposite direction} of the well-studied minimal controllability problems.  It is proven that all these problems are NP-hard, both for a numerical system and a structured system, {even restricted to some special cases.} {It is also shown that the first problems both for a numerical system and a structured one share a cardinality-constrained submodular minimization structure, for which there is no known constant or logarithmic factor approximation in general.} On the other hand, for the first two problems, under a reasonable assumption often met by practical systems, that the eigenvalue geometric multiplicities of the numerical systems or the matching deficiencies of the structured systems are bounded by a constant, {by levering the rank-one update property of the involved rank function}, it is possible to obtain the corresponding optimal solutions by {traversing a subset of the feasible solutions. We show such a method has polynomial time complexity in the system dimensions and the number of sensors.}

\end{abstract}
\begin{IEEEkeywords}
Observability robustness, secure estimation, structured system, computational complexity, {rank-one update} %,  submodular function minimization
\end{IEEEkeywords}

%\keywords{Controllability robustness, actuator removals, security analysis, recursive tree, submodular function minimization, complexity}

%%%%%%%%%%%%%%%%%%%%%%%%%%%%%%%%%%%%%%%%%%%%%%%%%%%%%%%%%%%%%%%%%%%%%%%%%%%%%%%%
\section{INTRODUCTION} \label{intro}
Modern control systems are often equipped with computer-based components, which could be vulnerable to cyber attacks \cite{wood2002denial}. As such, security has become a more and more important issue in control and estimation of cyber-physical systems, such as chemical processes, power grids, and transportation networks \cite{fawzi2014secure,shoukry2015event,mitra2019byzantine}. Typical attacks could be imposed on actuators, sensors or controllers of a control system. Among the related issues, one problem has attracted researchers' interest, which is that, under what condition it is possible to uniquely reconstruct the states of a system by observations of system outputs in the presence of adversarial sensor attacks \cite{chong2015observability,fawzi2014secure,shoukry2015event}.  Particularly, it is found in \cite{fawzi2014secure,chong2015observability} that, {\emph{unique state reconstruction of a linear time invariant (LTI) system under $s\in {\mathbb N}$ attacks (i.e., the number of sensors that are affected by attack signals) over a sufficiently long time horizon is possible, if and only {{if}} that system remains observable after the removal of arbitrary sets of $2s$ sensors ({{see \citep[Prop. 2]{fawzi2014secure} and \citep[Theo. 1]{chong2015observability} for details}}). }} Such condition is called \emph{$2s$-sparse observability} in \cite{shoukry2015event}, and has been extended to nonlinear systems and the distributed secure state estimation scenarios
\cite{kim2018detection,mitra2019byzantine}.  For a given system and $s\in \mathbb N$, verifying the aforementioned condition seems combinatorial.
This motivates our first problem on observability robustness under sensor failures.

%To be specific, consider the following linear time invariant (LTI) system under attack:
%\begin{equation}\label{ss} \dot x(t)=Ax(t), y(t)=Cx(t)+e(t),\end{equation}
%where $A\in {\mathbb R}^{n\times n}$, $C\in {\mathbb R}^{r\times n}$, $x(t)\in {\mathbb R}^{n},y(t)\in {\mathbb R}^r$ are respectively the state and output vectors, and $e(t)\in {\mathbb R}^r$ is the attack signal injected to sensors. The vector $e(t)$ is sparse in the sense that only the components corresponding to the attacked sensors can be nonzero, $t\ge 0$, but identities of the attacked sensors are not known.

More precisely, we will consider the following problem: {\emph{given an LTI system, determine the minimum number of sensors whose removal { can make the resulting system unobservable.}}}
We call this problem the {\emph{minimal sensor removal observability problem (minSRO)}}. Here, sensor removal could also result from denial-of-service attacks on a specific set of sensors \cite{wood2002denial}, which means that the attacked sensors cannot output any system information.

{ In addition to sensor removals, we will also briefly consider another type of sensor failures, that is, a subset of state variables are prevented from being directly measured by the existing sensors (in other words, these state variables are `blocked' from the existing sensors). Such a scenario may happen, for example, when an agent in a multi-agent system (or in a cooperative localization/tracking task) moves into a region where it loses all communications with the existing sensors \cite{ChenKalman2019,pal2010localization}. And the dual scenario is that, in a leader-follower system, a follower loses all its direct connections with the leaders. To measure observability robustness in such scenarios, we will consider the problem of {\emph{determining the minimum number of state variables that need to be blocked from the existing sensors to make the resulting system unobservable}}, which is referred to as the {\emph{minimal state variable blocking observability problem (minSBO)}}.}

%On the other hand, the dual of the minSRO, i.e., determining the minimum number of actuators (inputs) whose removal can destroy system controllability, can be seen as the {inverse counterpart} of the {\emph{minimal controllability problem (MCP)}}\footnote{{Note that there is no one-to-one correspondence between solutions to the two problems.}}, which is formulated as selecting the minimum number of inputs {\emph{from a given set of inputs}} to ensure system controllability \cite{T2016On}. The MCP is proven to be NP-hard, and can be approximated gracefully using simple greedy algorithms \cite{A.Ol2014Minimal, T2016On}. An extension of the MCP is to determine the sparsest input matrices such that the system remains controllable against a desired level of actuator removals \cite{Sergio_Pequito_2017_robust}. {The MCP in \cite{A.Ol2014Minimal,pequito2016minimum} could alternatively be formulated as selecting the minimum number of state variables to directly actuate by inputs to ensure (structural) controllability. In this sense, the dual of the minSBO could also be regarded as the inverse counterpart of the MCP in \cite{A.Ol2014Minimal,pequito2016minimum}.}

On the other hand, the dual of the minSRO, i.e., determining the minimum number of actuators (inputs) whose removal can destroy system controllability, is {in the opposite direction} of the well-studied {\emph{minimal controllability problem (MCP)}}, which is formulated as selecting the minimum number of inputs {\emph{from a given set of inputs}} to ensure system controllability \cite{T2016On}. The MCP is proven to be NP-hard and can be approximated gracefully using simple greedy algorithms \cite{A.Ol2014Minimal, T2016On}. An extension of the MCP is to determine the sparsest input matrices such that the system remains controllable against a desired level of actuator removals \cite{Sergio_Pequito_2017_robust}. {The MCP in \cite{A.Ol2014Minimal,pequito2016minimum} could alternatively be formulated as selecting the minimum number of state variables to be directly actuated to ensure (structural) controllability. In this sense, the dual of the minSBO is also opposite to the MCP in \cite{A.Ol2014Minimal,pequito2016minimum}.} Although the MCP has been extensively explored recently \cite{A.Ol2014Minimal, T2016On,Sergio_Pequito_2017_robust,pequito2016minimum}, little attention has been paid to its opposition direction.

{In the literature, structural controllability often serves as an alternative notion for controllability in the generic sense, owing to its nice property of relying only on the zero-nonzero patterns of system matrices rather than their exact values \cite{Y.Y.2011Controllability,generic}}. And because of this, many works have studied controllability/observability robustness under the structured system framework \cite{commault2008observability,Rahimian2013Structural,zhang2017edge,zhang2019minimal,lou2018toward,Y.Y.2011Controllability}. For example, \cite{commault2008observability} discusses observability preservation under sensor removals, and \cite{Rahimian2013Structural} studies controllability preservation under simultaneous failures in both the edges and nodes. Controllability robustness is also measured by the number of additional inputs needed for controllability under edge/node removals in \cite{Y.Y.2011Controllability,lou2018toward}. The above works mainly focus on the classifications of edges and nodes according to effects of their failures on system controllability/observability, rather than the optimization problems. A similar problem to the minSRO under structured system framework with {nonuniform} actuator costs has been proven to be NP-hard in \cite{zhang2019minimal}. However, no efficient algorithm was provided therein. {Another line to study controllability robustness is using
 strong structural controllability (SSC) \cite{mayeda1979strong}, as SSC could be regarded as
the ability to preserve system controllability for a structured system under arbitrary perturbations without changing the zero-nonzero patterns of the system matrices. In this line, \cite{mousavi2018structural} and \cite{jia2020unifying} respectively develop the conventional SSC to undirected networks and to allowing entries that can take arbitrary values. SCC preservation under structural perturbations is recently studied in \cite{mousavi2020strong}.}

In this paper, we consider the minSRO and the minSBO, as well as their structure counterparts, from a computational perspective.  The main contributions are as follows: 1) We prove that the minSROs are NP-hard both for a numerical system and a structured system, {even restricted to some special cases}. This confirms that verifying the $2s$-sparse observability condition in \cite{fawzi2014secure,shoukry2015event} is computationally intractable; {2) Following 1), we show the minSBOs are NP-hard both for a numerical system and a structured one; 3) We reveal that the minSROs both for numerical and structured systems share a cardinality-constrained submodular minimization structure that is hard to approximate in general; 4) Because of 3), instead of finding approximate algorithms, we turn to providing scenarios where polynomial-time exact algorithms\footnote{Exact algorithms are algorithms that return exactly the optimal solutions.} may exist. Under a reasonable constant bound assumption on the geometric multiplicities of system eigenvalues or matching deficiencies of system graphs, we show the minSROs can be solved in polynomial time by traversing a subset of the feasible solutions, both for numerical and structured systems.} Particularly, this means, in contrast to the MCP in \cite{A.Ol2014Minimal,T2016On}, the computational intractability of the minSRO for a numerical system is essentially caused by the increase of system eigenvalue geometric multiplicities, rather than that of numbers of system states or sensors.

{It is worth noting that independent of our work (see the preliminary version \cite{zhang2018arxiv}), \cite{Mao2019} recently shows the minSRO for a numerical system is NP-hard using a slightly different construction technique. One extension of the technique in our previous work \cite{zhang2019minimal} may also lead to the NP-hardness of the minSRO for a structured system. We will discuss their relations with this paper detailedly in the corresponding sections.}

The rest of this paper is organized as follows. Section II provides the problem formulations, and Section III presents some preliminaries. Sections IV and V give the computational complexity and algorithms for the considered problems, respectively. Section VI studies the structured counterparts of the minSRO and minSBO, with Section VII providing some numerical results. The last section concludes this paper.

% and its structure conterpart problem. It is proven that determining the minimal number of sensors whose removal destroy system observability is generally NP-hard. Nevertheless, a pseudo-polynomial time algorithm is provided for solving this problem on systems with bounded maximum geometric multiplicities. This algorithm uses traversals over a recursive tree built from the eigenspaces of the system state transition matrix and the output matrix. This result indicates that, in contrast to the minimal controllability problems, the computational intractability of the minSRO problem is essentially caused by the increase of maximum geometric multiplicities of system matrices, rather than that of the system dimensions or number of sensors.  Moreover, this algorithm is extended to the case where each sensor has a non-negative cost. Considering that this algorithm requires accurate eigenspace decomposition of the system matrices, we then focus on a similar problem on structured systems or generic systems, i.e., systems which can be represented by graphs and the nonzero entries of whose system matrices are unknown but can take values independently. It is proven that such problem is still NP-hard for the structured systems. Nevertheless, a pseudo-polynomial time algorithm is provided on a structured system with bounded matching deficiencies.

%Symbols ${\mathbb R}, {\mathbb C}$ and $\mathbb N$ denote the sets of real, complex and natural numbers, respectively.
%{\emph{Related Work:}}

{\emph{\bf { Notations}:}} For an $m\times n$ ({numerical or structured}) matrix $M$, upon letting $J_1\subseteq\{1,...,m\}$ and $J_2\subseteq \{1,...,n\}$, $M_{J_1,J_2}$ denotes the submatrix of $M$ whose rows are indexed by $J_1$ and columns are indexed by $J_2$, $M_{J_1}$ the submatrix consisting of {\emph{rows}} indexed by $J_1$, and $M_{[J_2]}$ the submatrix consisting of {\emph{columns}} indexed by $J_2$. By ${\bf col}\{X_i|_{i=1}^N\}$ (resp. ${\bf diag}\{X_i|_{i=1}^N\}$) we denote the composite (resp. diagonal) matrix with its $i$th row block (resp. diagonal block) being $X_i$, and $\overline{{\rm span}}(X)$ the space spanned by {\emph{row}} vectors of $X$. The combinatorial number ${\tiny{\left(\begin{array}{l} n\\m \end{array}\right)}}$ is also denoted by ${\bf C}_n^{m}$ ($m\le n$).

\vspace{-0.4em}
\section{Problem Formulations}
%\vspace{-0.5em}
Consider the following LTI system
\begin{equation}
\label{plant Eq}
\dot{x}(t)=Ax(t)+Bu(t),\\~
y(t)=Cx(t)
\end{equation}
where $x(t) \in \mathbb{R}^n$, $u(t) \in \mathbb{R}^m$ and $y(t)\in \mathbb{R}^{r}$ are respectively the state, input and output vectors, $A \in \mathbb{R}^{n\times n}$, $B=[b_1,...,b_m] \in \mathbb{R}^{n \times m}$  and $C={{\bf col}\{c_i|_{i=1}^r\}}\in \mathbb{R}^{r\times n}$ are respectively the state transition, input and output matrices, with $b_i\in {\mathbb R}^{n}$, $c_i\in {\mathbb R}^{1\times n}$. Without loss of generality, {\emph{assume that every row of $C$ is nonzero}}. Throughout this paper, $V\doteq\{1,...,r\}$. %  and $C_S$ denotes the submatrix of $C$ consisting of rows of $C$ indexed by $S\subseteq V$.

%System (\ref{plant Eq}) is said to be controllable, if for any two states $x_0, x_1 \in{{\mathbb{R}}^n}$, there exists an input $u(t)$ that can drive the system states from $x_0$ to $x_1$ in finite time. We just simply say $(A,B)$ is controllable if System (\ref{plant Eq}) is controllable.

%In the field of cyber-physical security, it is of both theoretical and practical interest to know whether a given system like (\ref{plant Eq}) can preserve its observability under the removal of a cardinality-constrained set of sensors \cite{fawzi2014secure,shoukry2015event,mitra2019byzantine, an2019distributed}. As mentioned earlier, this is not only related to the ability of unique reconstruction of system states under adversarial sensor attacks, but also significant to measure observability robustness/resilience under denial-of-service attacks on sensors.
%the maximum number of adversaries that can be tolerated.
%For practical cyber-physical systems, we may want to know whether a given system like (\ref{plant Eq}) can preserve its controllability under denial of service attacks on a cardinality-constrained set of sensors.
We say System (\ref{plant Eq}) is $s$-robust observable, $s\in \mathbb{N}$, if under the removal of arbitrary $s$ sensors, the resulting system is still observable.
To measure such observability robustness/resilience, the following minSRO is considered.
\[ \begin{array}{*{20}{l}}
{\mathop {\min }\limits_{S \subseteq V} {\kern 1pt} {\kern 1pt} {\kern 1pt} {\kern 1pt} {\kern 1pt} \left| S \right|}\tag{Problem 1}  \\
{{\rm{s.t. }}\ (A,{C_{V\backslash S}}){\kern 1pt} {\kern 1pt} {\kern 1pt}  {\kern 1pt} {\kern 1pt} {\kern 1pt} {\rm{unobservable}}}
\end{array}\]
Denote the cardinality of the optimal solution to Problem 1 by $r_{\min}$. Then, System (\ref{plant Eq}) is $s$-robust observable for any integer $s\le r_{\min}-1$. Moreover, according to \cite{chong2015observability,fawzi2014secure}, it is possible to uniquely reconstruct the system states under $s'$ sensor attacks with $s'\le \left\lfloor {\frac{r_{\min}-1}{2}} \right\rfloor$, where $\left\lfloor {\bullet} \right\rfloor$ denotes the floor of a scalar.

Next, we will consider the structure counterpart of Problem 1. This is motivated by the fact that the exact values of system matrices may sometimes be hard to know in practice due to modeling errors or parameter uncertainties. Instead, the sparsity patterns of system matrices are often easily accessible. Under this circumstance, one may be interested in the generic properties of a system, i.e., properties that hold almost everywhere in the corresponding parameter space \cite{generic}. Controllability and observability are such generic properties.

%It is known that, if a structured system is structurally observable, then almost all of its numerical realizations correspond to observable systems; otherwise, none of its system realizations is observable.

To be specific, a {\emph{structured matrix}} is a matrix whose entries are either fixed zero or free parameters. Denote the set of all $n\times m$ structured matrices by $\{0,*\}^{n\times m}$, where $*$ denotes the free parameters. A {\emph{structured system}} is a system whose system matrices are represented by structured matrices. This system is said to be {\emph{structurally observable}}, if there exists one numerical realization of it (i.e., a pair of numerical matrices that has the same sparsity pattern as the structured matrices) that is observable. Let $\bar A$ and $\bar C$ be structured matrices consistent with the sparsity patterns of $A$ and $C$. We consider the following problem.
\[\begin{array}{*{20}{l}}
{\mathop {\min }\limits_{S \subseteq V} {\kern 1pt} {\kern 1pt} {\kern 1pt} {\kern 1pt} {\kern 1pt} \left| S \right|}\tag{Problem 2} \\
{{\rm{s.t.}}\ (\bar A,{\bar C_{V\backslash S}}){\kern 1pt} {\kern 1pt} {\kern 1pt} {\kern 1pt} {\kern 1pt} {\kern 1pt} {\rm{structurally {\kern 2pt} unobservable}}}
\end{array}\]

%\[\begin{array}{*{20}{l}} \begin{aligned}
%{Problem \ 2:}\ &{\mathop {\min }\limits_{S \subseteq V} {\kern 1pt} {\kern 1pt} {\kern 1pt} {\kern 1pt} {\kern 1pt} \left| S \right|} \\
%&{{\rm{s.t.}}\ (\bar A,{\bar C_{V\backslash S}}){\kern 1pt} {\kern 1pt} {\kern 1pt} {\kern 1pt} {\kern 1pt} {\kern 1pt} {\rm{structurally {\kern 2pt} unobservable}}}\end{aligned}
%\end{array}\]

By definition of structural observability, the optimal value of Problem~2 is an {\emph{upper bound}} of that of Problem 1.
As observability is a generic property, it is expected that the optimal values of these two problems tend to coincide with each other in randomly generated systems, which will also be demonstrated by the numerical experiments in Section \ref{numer_example}.

% Our numerical experiments in Section xx will show that, the optimal values of these two problems coincide with each other in most randomly generated systems.

%[As a by-product of our derivations for Problems 1 and 2]
{                                                                                                                                                                                                                                                                                                                                                                                                                                                                                                                                                                                                                                                                                                                                                                                                                                                                                                                                                                We will also briefly consider another type of sensor failures, that is, a subset of state variables are blocked from the existing sensors (corresponding to that the associated nonzero columns of $C$ are transformed into zero).  Let $\hat V \subseteq \{1,...,n\}$ be the set of indices of nonzero columns of $C$ (resp. $\bar C$). For $S\subseteq \hat V$, denote $C^{S}$ (resp. $\bar C^{S}$) by the matrix obtained from $C$ (resp. $\bar C$) by preserving its columns indexed by $S$ and zeroing out the rest. To measure such observability robustness, the minSBO and its structured counterpart can be formulated respectively as
%\label{problem 3add}
\[\begin{array}{*{20}{l}}
{\mathop {\min }\limits_{S \subseteq \hat V} {\kern 1pt} {\kern 1pt} {\kern 1pt} {\kern 1pt} {\kern 1pt} \left| S \right|} \tag{Problem 3}\\
{{\rm{s.t.}}\ (A,{C^{\hat V \backslash S}}){\kern 1pt} {\kern 1pt} {\kern 1pt}  {\kern 1pt} {\kern 1pt} {\kern 1pt} {\rm{unobservable}}}
\end{array}\]and% \label{problem 4add}
\[\begin{array}{*{20}{l}}
{\mathop {\min }\limits_{S \subseteq \hat V} {\kern 1pt} {\kern 1pt} {\kern 1pt} {\kern 1pt} {\kern 1pt} \left| S \right|}\tag{Problem 4} \\
{{\rm{s.t.}}\ (\bar A,{\bar C^{\hat V \backslash S}}){\kern 1pt} {\kern 1pt} {\kern 1pt} {\kern 1pt} {\kern 1pt} {\kern 1pt} {\rm{structurally \ unobservable}}}
\end{array}\]
Note that when $S=\hat V$, the resulting system will become (structurally) uncontrollable, which justifies that Problems 3 and 4 are always feasible.
As mentioned earlier, Problems 3 and 4 can measure observability robustness in a cooperative localization/tracking system where partial individuals lose all communications with the sensors \cite{ChenKalman2019,pal2010localization}.
}

%Furthermore, a variant of Problem 1 is to determine the minimum dimension of observable subspaces of a given system under cardinality-constrained actuator failures.  This problem asks for the worst damage in observable subspaces that the failures of a set of sensors with cardinality upper bound can cause. Denote the observability matrix of a pair $(A,B)$ by ${\mathcal{C}}(A,B)$, i.e., ${\mathcal{C}}(A,B)=[B,AB,...,A^{n-1}B]$. Then, this problem can be formulated as
%
%{\begin{problem} \label{problem 2}
%Given $(A,B)$ in (\ref{plant Eq}) and $l\in \mathbb{N}$
%\[ {\mathop {\min }\limits_{S \subseteq V, |S|\le l} {\kern 2pt}  {\rm rank} {\kern 1pt} {\mathcal{C}}(A,B_{V\backslash S})} \]
%\end{problem}

% {{The main purpose of this paper is to reveal the complexity of Problem 1,  develop algorithms for its related problems, and apply them to security analysis of attacks on sensors.}}

\section{Preliminaries} \label{sec_preliminary}
This section presents some necessary preliminaries on (structural) observability and graph theories.

%The following well-known PBH test gives a necessary and sufficient condition for System (\ref{plant Eq}) to be observable.
{\begin{lemma} \label{lemma 1}[PBH test,\cite{Kailath_1980}]
System (\ref{plant Eq}) is observable, if and only if for each $\lambda\in \mathbb{C}$,  there exists no $x\in {\mathbb{C}}^n$, $x\ne 0$, such that $Ax=\lambda x$ and $Cx=0$.
\end{lemma}}

Given System (\ref{plant Eq}), suppose $A$ has $p\le n$ distinct eigenvalues, and denote the $i$th one by $\lambda_i$. Let $k_i$ be the dimension of the null space of $\lambda_i I- A$, i.e., $k_i$ is the geometric multiplicity of $\lambda_i$. Let $X_i$ be an eigenbasis of $A$ associated with $\lambda_i$, that is, $X_i$ consists of $k_i$ vectors which are linearly independent spanning the null space of $\lambda_i I- A$. With these definitions, the following lemma is immediate from the PBH test.

\begin{lemma}[\cite{Kailath_1980}] \label{corollary 1}
System (\ref{plant Eq}) is observable, if and only if $CX_i$ is of full column rank for $i=1,...,p$.
\end{lemma}

%Unlike the PBH test for observability which relies on eigenspace decomposition of the system matrices,
Criteria for structural observability often have explicit graphical presentations.  To this end, for System (\ref{plant Eq}), let $\mathcal{X}=\{x_1,...,x_n\}$, $\mathcal{Y}=\{y_1,...,y_m\}$ denote the sets of state vertices and output vertices respectively. Denote the edges by $\mathcal{E}_{\mathcal{X},\mathcal{X}}=\{(x_i,x_j): \bar A_{ji}\ne 0\}$,
$\mathcal{E}_{\mathcal{X},\mathcal{Y}}=\{(x_j,y_i): \bar C_{ij}\ne 0\}$.  Let $\mathcal{D}(\bar A,\bar C)= (\mathcal{X}\cup \mathcal{Y}, \mathcal{E}_{\mathcal{X},\mathcal{X}}\cup \mathcal{E}_{\mathcal{X},\mathcal{Y}})$ be the system digraph associated with $(\bar A, \bar C)$, and $\mathcal{D}(\bar A)=(\mathcal{X},\mathcal{E}_{\mathcal{X},\mathcal{X}})$. A state vertex $x\in\mathcal{X}$ is said to be output-reachable, if there exists a path from $x$ to an output vertex $y\in\mathcal{Y}$ in $\mathcal{D}(\bar A,\bar C)$.

The generic rank of a structured matrix $M$, denoted by ${\rm grank}(M)$, is the maximum rank $M$ can achieve as the function of its free parameters.
The bipartite graph associated with a  matrix $M$ is given by ${\mathcal B}(M)=({\mathcal R}(M),{\mathcal C}(M),{\mathcal E}(M))$, where the left vertex set ${\mathcal R}(M)$ (resp. right vertex set ${\mathcal C}(M)$) corresponds to the row (resp. column) index set of $M$, and the edge set corresponds to the set of nonzero entries of $M$, i.e., ${\mathcal E}(M)=\{(i,j): i\in {\mathcal R}(M), j\in {\mathcal C}(M), M_{ij}\ne 0\}$. A matching of a bipartite graph is a subset of its edges among which any two do not share a common vertex. A vertex is said to be unmatched associated with a matching if it is not in edges of this matching. The maximum matching is the matching with the largest number of edges among all possible matchings. It is known that ${\rm grank}(M)$ equals the cardinality of the maximum matching of ${\mathcal B}(M)$ \cite{Murota_Book}.

\begin{lemma}[\cite{generic}]\label{Lemma 1}
The pair $(\bar A, \bar C)$ is structurally observable, if and only if  (a) every state vertex is output-reachable in $\mathcal{D}(\bar A,\bar C)$, and
             (b) ${\rm grank}([\bar A^{\intercal}, \bar C^{\intercal}])=n$.
\end{lemma}

We shall call (a) of Lemma \ref{Lemma 1} the {\emph{output-reachability condition}}, and (b) the {\emph{matching condition}}.

\section{Computational Complexity of Problems 1 and 3}
A natural question is whether Problem 1 is solvable in polynomial time. In this section, we will give a negative answer, {following which the NP-hardness of Problem 3 is immediate.}

\begin{theorem}\label{theorem 1}
Problem 1 is NP-hard, {even with dedicated sensors (i.e., each sensor measures only one state variable)}.
\end{theorem}

Our proof is based on the reduction from the {\emph{linear degeneracy problem}}, defined in Definition \ref{def-linear}. In \cite{Khachiyan_1995}, it is proven that this problem is NP-complete, and there are infinitely many { integer} matrices associated with which this problem is NP-complete.

\begin{definition}[\cite{Khachiyan_1995}] \label{def-linear}
The linear degeneracy problem is to determine whether a given $d\times p$ rational matrix $F=[f_{1},...,f_{p}]$ contains a degenerate submatrix of order $d$, i.e., $\det [f_{i_1},...,f_{i_d}]=0$, for some $\{i_1,...,i_d\} \subseteq \{1,...,p\}$.
\end{definition}

{\bf{Proof of Theorem \ref{theorem 1}:}} Let $X$ be an arbitrary $k\times n$ { integer} matrix ($k< n$) with full row rank, and $\alpha_1=\max \{|X_{ij}|\}$.  %Notice that the full row rank constraint does not alter the NP-completeness of the linear degeneracy problem associated with $X$.
Let $X^\bot$ be the basis matrix of the null space of $X$, i.e., $X^{\bot}$ is an $n\times (n-k)$ matrix spanning the null space of $X$. $X^{\bot}$ can be constructed via the Gaussian elimination method in polynomial time ${O}(n^3)$ \cite{Geor1993Graph}. Moreover, let $\alpha_2=\max \{|X^{\bot}_{ij}|\}$, and $\alpha_{\max}=\max\{\alpha_1, \alpha_2 \}$.  Notice that the entries of $X^{\bot}$ are rational, and $X^{\bot}$ multiplied by any nonzero scalar is still a basis matrix of the null space of $X$. Hence, the encoding length of $\alpha_{2}$ (i.e., $\log_2\alpha_{2}$) can be polynomially bounded by $k$ and $n$; so is $\alpha_{\max}$. Define $H(\eta)$ as
$$H(\eta ) = \left[X^{\intercal},X^{\bot}+\eta{\bf 1}_{n\times (n-k)}\right],$$
where ${\bf 1}_{n\times (n-k)}$ denotes the $n\times (n-k)$ matrix whose entries are all one. Then, clearly, $\det H(0)\ne 0$.
We are to show that one can {\emph{deterministically}} find a scalar $\eta ^*$ satisfying $\eta ^* > \alpha_{2}$ and $\det H(\eta ^*)\ne 0$ in polynomial time. Notice that $\det H(\eta)$ is a polynomial of $\eta$ with degree at most one as the coefficient matrix of $\eta$ in $H(\eta)$ is rank-one. Hence, $\det H(\eta)$ can be expressed as $\det H(\eta)=[\det H(1)-\det H(0)]\eta+\det H(0)$.\footnote{{This can be seen as follows. Note $\det H(\eta)$ can be expressed as $\det H(\eta)=\beta_0+\beta_1\eta$ ($\beta_0, \beta_1\in {\mathbb C}$). Upon solving the system of equations $\det H(0)=\beta_0$ and $\det H(1)=\beta_0+\beta_1$, we get $\beta_0=\det H(0)$ and $\beta_1=\det H(1)-\det H(0)$, which results in the proposed expression.}} {Note also $\det H(\eta)$ has at most one root for $\eta$ as $\det H(0)\ne 0$.}  %Hence, $\det H(\eta)$ can be expressed as $\det H(\eta)=\beta_0+\beta_1\eta$ ($\beta_0, \beta_1\in {\mathbb C}$). Upon substituting $\det H(0)$ and $\det H(1)$ into that expression, we get $\beta_0=\det H(0)$ and $\beta_1=\det H(1)-\det H(0)$, which results in} $\det H(\eta)=[\det H(1)-\det H(0)]\eta+\det H(0)$.
Therefore, in arbitrary set consisting of $2$ distinct rational numbers which are bigger than $\alpha_{2}$, there must exist one $\eta^*$, such that $\det H(0) (1-\eta^*)+\det H(1)\eta^*\ne 0,$
leading to $\det H(\eta ^*)\ne 0$. %Furthermore, dy the definition of determinant, it can be verified that both $\det H(0)$ and $\det H(1)$ have encoding lengths (which are  $n\log_2\alpha_{\max} +n \log_2n$ and $n\log_2(\alpha_{\max}+1) +n \log_2n$ respectively) polynomially bounded by $n$.
Moreover, dy the definition of determinant,  it holds $\det H(0)\le \alpha_{\max}^nn^n$ and $\det H(1)\le (\alpha_{\max}+1)^nn^n$, which means $\det H(0)$ and $\det H(1)$ have encoding lengths (being $n\log_2\alpha_{\max} +n \log_2n$ and $n\log_2(\alpha_{\max}+1) +n \log_2n$ respectively) polynomially bounded by $n$.
%After determining an $\eta ^*$ satisfying $\eta ^* > \alpha_{2}$ and $\det H(\eta ^*)\ne 0$,
Afterwards, let matrices $P=H(\eta^*)$, $\Gamma={\bf{diag}}\{I_k,2,3,...,n-k+1\}$, and construct the system~as
$$A=P\Gamma P^{-1}, C=I_n.$$
Since the encoding lengths of entries of $P$ are polynomially bounded by $n$, its inversion $P^{-1}$ can be computed in polynomial time and has polynomially bounded encoding lengths too. Hence, $A$ can be computed in polynomial time.

We claim that the optimal value of Problem 1 associated with $(A, C)$ is no more than $n-k$, if and only if there exists an $k\times k$ submatrix of $X$ which has zero determinant.

Indeed, from the construction of $A$, the $(k+i)$th column of $P$, denoted by $P_{[k+i]}$, is the eigenvector of $A$ associated with the eigenvalue $i+1$, $i=1,...,n-k$. Notice that, all entries of $P_{[k+i]}$ are nonzero as $\eta^*>\alpha_2$. Hence, all $n$ rows of $C$ need to be removed, such that the resulting $C_{V\backslash S}P_{[k+i]}$ fails to be of full column rank, where $S\subseteq V$. From Lemma \ref{corollary 1}, the optimal value of Problem 1 associated with $(A,C)$ then equals the minimum number of rows
(resp., columns) whose removal from $X^{\intercal}$ ({resp.}, $X$) makes the resulting matrix fail to have full column (row) rank. If such value is no more than $n-k$,  there must exist an $k\times k^{+}$ submatrix of $X$ which is not of full row rank for some $k^+\ge k$. Then, the linear degeneracy problem associated with $X$ is yes.

Conversely, suppose there is an $k\times k$ submatrix of $X$ with zero determinant, and denote it by $X_{[\bar S]}$, $\bar S \subseteq V$. Then, one just needs to remove the columns indexed by $V\backslash {\bar S}$ from $X$, such that the resulting $X_{[\bar S]}$ fails to be of full row rank. Hence, the optimal value of Problem 1 is no more than $|V\backslash {\bar S}|=n-k$. The fact that the linear degeneracy problem associated with $X$ is NP-complete yields the NP-hardness of Problem~1. $\hfill\blacksquare$

The above result indicates that, there is in general no efficient way to verify whether a given system is $s$-robust observable for a given $s\in {\mathbb N}$ unless $P=NP$, {even when $C=I_n$}. This confirms the conjecture implicitly made in \cite{chong2015observability,fawzi2014secure,shoukry2015event} where they tend to think that computing the maximal number of sensor attacks that can be tolerated is computationally intractable (see Section \ref{intro}). {It is worth noting that an independent work \cite{Mao2019} from the present paper (see \citep[Coro. 3]{zhang2018arxiv}) also shows it is NP-hard to verify the $s$-robust observability condition, which uses a similar technique but without the dedicated output configuration constraint. The dedicated output/input configuration might be common both in theoretical studies and practical implementations \cite{T2016On,Sergio_Pequito_2017_robust,A.Ol2014Minimal}. Therefore,  it would be desirable to make clear the complexity of Problem 1 even in this case.} Furthermore, with $C=I_n$, removing one sensor corresponds to that {\emph{exactly}} one state variable loses its direct measurements. Hence, Theorem \ref{theorem 1} readily leads to the following corollary.

{\begin{corollary} { Problem 3 is NP-hard.}
\end{corollary}} %Given a system $(A,C)$, it is NP-hard to determine the minimal number of state variables that need to be prevented from being directly measured by the existing sensors (i.e., the number of nonzero columns of $C$ that need to be zeroed out), such that the resulting system becomes unobservable. % the minimal number of state variables that need to be blocked their direct accesses to the existing sensors % need to be prevented from being directly measured by the existing sensors

Next,  by the duality between controllability and observability, it is easy to obtain the following corollary on controllability robustness under actuator failures/removals.

{\begin{corollary} For System (\ref{plant Eq}), it is NP-hard to determine the minimum number of actuators whose removal can make the resulting system uncontrollable.
\end{corollary}}

Taking the input and output into consideration together, we have the following conclusion.
\begin{corollary} For System (\ref{plant Eq}), it is NP-hard to determine the minimal total number of sensors and actuators whose removal makes the resulting system neither controllable nor observable.
\end{corollary}

{\bf{Proof:}} Construct $(A,C)$ as suggested in the proof of Theorem \ref{theorem 1}, i.e., $A=P\Gamma P^{-1}$, $C=I_n$. Construct $B$ such that $B\in \mathbb{R}^{n\times k}$ and $(A,B)$ is controllable. Since the eigenbases $P$ of $A$ are available, the matrix $B$ can be constructed in polynomial time as suggested in \cite{Tong2017Minimal}. Notice that $k$ is the maximum geometric multiplicity of eigenvalues of $A$. According to \cite{Tong2017Minimal}, the minimal number of inputs (actuators) that ensures controllability of the associated system equals $k$. Hence, removing any one of these $k$ actuators can make the system uncontrollable. As a result, the minimal total number of actuators and sensors whose removal causes uncontrollability and unobservability equals $r_{\min}+1$, where $r_{\min}$ is the optimum of the minSRO with $(A,C)$. The latter problem is shown to be NP-hard in Theorem \ref{theorem 1}. The required result follows. $\hfill\blacksquare$

\section{Structure and Algorithm for Problem 1}
%In this section, we first give a polynomial time algorithm for Problem 1 on systems with bounded eigenvalue geometric multiplicities. Then, we extend this algorithm to the case where each sensor has a non-negative cost. We shall assume that a collection of eigenbases of $A$ are computationally available, and denote them by $X_i|_{i=1}^p$. The symbols $p$, $k_i|_{i=1}^p$ are defined in Section~III.

%[In the above two sections, we have demonstrated the intractability of Problems 1-4, even restricted to some special cases. In this section, we deal with the algorithmic perspectives. We show that Problems 1 and 3 share a cardinality-constrained submodular minimization (CCSM, i.e., minimizing a submodular function subject to an upper/lower bound on cardinality) structure, which indicates they might be hard to approximate in general. On the other hand, under a
%reasonable assumption often met by practical engineering systems (see Remark x,x), we show that Problems 2-4 can be solved in polynomial time using a similar framework.]

In this section, we first show that Problem 1 has a cardinality-constrained submodular minimization (CCSM, i.e., minimizing a submodular function subject to an upper/lower bound on cardinality) structure, which indicates it might be hard to approximate. We then show Problem 1 can be solved in polynomial time under a constant bound assumption on the system eigenvalue geometric multiplicities.

We briefly introduce the notion of submodularity here; see \cite{Murota_Book} for more details. Let $V$ denote a finite set. A function $f: 2^{V}\rightarrow {\mathbb R}$ that takes a subset of $V$ as the input and outputs a real value is called submodular, if for any $S,T\subseteq V$, $f(S)+f(T)\ge f(S\cap T)+f(S\cup T)$.
\vspace*{-3mm}
\subsection{The CCSM Structure of Problem 1} \label{CCSM-sec-pro1}
Let ${\cal O}(A,C)$ be the observability matrix of $(A,C)$, i.e., ${\cal O}(A,C)={\bf col}\{C,CA,...,CA^{n-1}\}$. Then, it can be seen that Problem 1 can be solved by invoking the following problem for at most $O(\log_2r)$ times using bisection on $l\in \{0,...,r\}$ until the corresponding objective is less than $n$
\begin{equation}\label{ccsm1} {\mathop {\min }\limits_{S \subseteq V, |S|\le l} {\kern 2pt}  {\rm rank} {\kern 1pt} {\mathcal{O}}(A,C_{V\backslash S})}.\end{equation}
Define $f(S):2^V\rightarrow \mathbb{N}$ as $f(S)={\rm rank} {\kern 1pt} {\mathcal{C}}(A,C_{S})$. It is known that $f(S)$ is submodular on $S\subseteq V$ \citep[Theo. 7]{T2016On}. From the property of submodularity \cite{Murota_Book},  $f(V\setminus S)$ is also submodular on S$\subseteq V$. This means that (\ref{ccsm1}) can be formalized as a CCSM problem. However, the general CCSM is hard to approximate, and there is no known constant or logarithmic factor for it. Specifically, \cite{2008Submodular} has shown that approximating the submodular function minimization over an upper bound on cardinality within a factor of $O(\sqrt{\frac{r}{{\rm ln} r}})$  is generally impossible, where $r\doteq |V|$. \footnote{We can also reformulate Problem 1 as minimizing $|S|$ subject to ${\rm rank} {\kern 1pt} {\mathcal{C}}(A,C_{V\backslash S})\le n-1$. From the bicriteria results in \cite{2008Submodular}, approximating the cardinality minimization subject to a monotone submodular function value upper bound better than $O(\sqrt{\frac{ln r}{r}})$ is generally impossible, that is, upon letting $r_{\rm {alg}}$, $r_{\min}$ and $r$ be respectively the returned value by a polynomial algorithm, the optimal value, and the cardinality of $V$, making $r-r_{\rm {alg}}$ closer to $r-r_{\min}$ than a factor of $O(\sqrt{\frac{ln r}{r}})$ is impossible.} Hence, in the following, we turn to providing a scenario where the
polynomial-time (in $n$ and $r$) exact algorithm exists for Problem 1.

%\subsection{Simple Dynamics Case}
%At the beginning, we consider the simple spectrum case where $A$ has no repeated eigenvalues.
%In this case, $\{X_i|_{i=1}^n\}$ becomes a collection of eigenvectors of $A$. For the $i$th eigenvalue of $A$, define $F_i$ as the collection of rows of $C$ that are not orthogonal to $X_i$. Since $X_i$ is a vector, $c_jX_i$ becomes a scalar, $j\in \{1,...,r\}$. One has
%$F_i=\{j:c_jX_i\ne 0,j\in \{1,...,r\}\}.$ Then, from Lemma \ref{corollary 1}, it is obvious that $$r_{\min}=\min\limits_{i\in\{1,...,n\}} |F_i|.$$
%Obviously, finding $r_{\min}$ can be done in polynomial time.  The above analysis may seem trivial. In the following we will extend it to the general case where the geometric multiplicities of eigenvalues of $A$ can be greater than one.
%\vspace*{-2mm}
\subsection{Algorithm with Bounded Geometric Multiplicities}
We shall assume that a collection of eigenbases of $A$ are computationally available. The symbols $p$, $X_i|_{i=1}^p$ and $k_i|_{i=1}^p$ are defined in Section~III. Define  $Y^{(i)}\doteq CX_i$, $i=1,...,p$.

We further assume that the geometric multiplicities of eigenvalues of $A$ are bounded by some fixed constant $\bar k\in \mathbb{N}$. That is, $k_i\le \bar k$, $i=1,...,p$, as $n$ and $r$ increase. For most practical engineering systems, this assumption is reasonable. In fact, it is found in \cite{Tao2017random} that, with probability tending to $1$, the symmetric matrices with random entries have no repeated eigenvalues, and Erdos-Renyi (ER) random graph has
simple spectrum asymptotically almost surely. { Besides, in the areas of multi-agent systems or distributed estimation, Laplacian matrices of undirected graphs with bounded\footnote{Hereafter, if not specified, by being bounded we mean being bounded by a (known) constant as the input size grows.} connected components, or of weakly-connected directed leader-follower networks with bounded number of leaders, tend to have the maximum geometric multiplicities of their eigenvalues that are more than one but bounded by the corresponding constants \cite{mesbahi2010graph}. Similar scenarios may arise in modeling composite/layer networks in which the system matrices are often (Cartesian) products of smaller factor graphs \cite{chapman2014controllability}.}

Let us focus on an individual eigenvalue $\lambda_i$, $i\in\{1,...,p\}$. Let $r_i$ be the minimum number of rows whose removal from
 $Y^{(i)}$ makes the remaining matrix fail to be full of column rank. To determine $r_i$, an exhaustive combinatorial search needs to compute the ranks of at most ${\bf C}_r^1+\cdots+{\bf C}_r^{r-k_i}\to O(2^r)$ submatrices, which increases exponentially with $r$ even when $k_i$ is bounded. Hence, the direct combinatorial search is not computationally efficient.

%{\footnotesize{$\left( \begin{array}{l}
%r\\
%1
%\end{array} \right) + \left( \begin{array}{l}
%r\\
%2
%\end{array} \right){\rm{ + }} \cdots {\rm{ + }}\left( {\begin{array}{*{20}{c}}
%r\\
%{r - {k_i }}
%\end{array}} \right) \to O({2^r})$}}

\begin{algorithm} % 算法开始
 {\small {{{{
\caption{: Algorithm for Problem 1} %算法的题目
\label{alg1} %算法的标签
\begin{algorithmic}[1] %此处的[1]控制一下算法中的每句前面都有标号
\REQUIRE Parameters $(A,C)$ of a controllable system  %, a permutation $g=\{k_1,...,k_p\}$ of the integer set $\{1,...,p\}$
  %输入条件(此处的REQUIRE默认关键字为Require，在上面已自定义为Input)
\ENSURE  The optimal solution $F_{\min}$ to Problem 1  %输出结果(此处的ENSURE 默认关键字为Ensure在上面已自定义为Output)
\STATE Calculate the eigenbases $\{X_i|_{i=1}^p\}$ of $A$.
\FOR{$\kappa =1$ to $p$}
\STATE  Initialize $\tau {\rm{ = }}0$, $T_1^{0} = \emptyset$, $L_{-1}=1$.
\WHILE{$\tau  < {k_\kappa }$}
%\COMMENT{// We use $\{c^{\tau}_0,c^{\tau}_1,...,c^{\tau}_{\rm end}\}$ to denote a sequence of sorted integers with updated lengths, where $c^{\tau}_{\rm end}$ denotes the last element of this sequence.}
\FOR{$i=1$ to $L_{\tau-1}$}
%\COMMENT{We use $\{c^{\tau}_0,c^{\tau}_1,...,c^{\tau}_{\rm end}\}$ to denote a sequence of sorted integers with updated lengths, where $c^{\tau}_{\rm end}$ denotes the last element of this sequence.}
%\STATE $\bar T_i^{\tau} =   T_i^{\tau} \cup \left\{{j:\overline{{\rm span}}(Y_{\{j\}}^{(\kappa )}) \subseteq \overline{{\rm span}}(Y_{T_i^{\tau}}^{(\kappa )}),j \in V\backslash T_i^{\tau}} \right\}.$
\STATE $\bar T_i^{\tau} =   T_i^{\tau} \cup \left\{{j:{\rm rank}(Y_{T_i^\tau\cup \{j\}}^{(\kappa )})={\rm rank}(Y_{T_i^{\tau}}^{(\kappa )}),j \in V\backslash T_i^{\tau}} \right\}.$
 \STATE    $\Omega _i^{\tau} = V\backslash \bar T_i^{\tau}$. Let $c_i^{\tau} = \sum\nolimits_{l = 1}^i {{\rm{|}}\Omega _l^{\tau}{\rm{|}}}$, $c_0^{\tau}{\rm{ = }}0$.
 \STATE Rewrite $\Omega _i^{\tau} \doteq \{ {j_1},...,{j_{|\Omega _i^{\tau}|}}\}$ and let $T_{c_{i - 1}^{\tau} + q}^{\tau  + 1}{\rm{ = }}\bar T_i^{\tau}\bigcup {{\rm{\{ }}{j_q}{\rm{\} }}}$ for $q = 1,...,|\Omega _i^{\tau}|$.
\ENDFOR
\STATE  $L_{\tau}=c^{\tau}_{L_{\tau-1}}$; $\tau=\tau +1$.
%\COMMENT{THIS IS}
\ENDWHILE
\STATE $T_{\max }^{[\kappa ]} = \arg \mathop {\max } \left\{ |\bar T_1^{{k_\kappa } - 1}|, \cdots ,|\bar T_{L_{{k_\kappa } - 2}}^{{k_\kappa } - 1}|\right\}$,\\~ ${F_\kappa } = V\backslash T_{\max }^{[\kappa ]}$.
 \ENDFOR%$ associated with $(V(\mathcal{T})\cup V_L(\mathcal{M})\cup V_R(\mathcal{M}), E(\mathcal{T})\cup \mathcal{M})$.
 \STATE Return $F_{\min}=\arg \mathop {\min} \left\{ |{F_1}|,...,|{F_p}|\right\}$ and ${r_{\min }} = |F_{\min}|$.
\end{algorithmic}}}
}}}
%{\footnotemark[1]{\emph{ Use $\{c^{\tau}_0,c^{\tau}_1,...,c^{\tau}_{\rm end}\}\subseteq \{1,...,r^{\tau}\}$ to denote a sequence of sorted integers with updated cardinality depending on $\tau$, where $c^{\tau}_{\rm end}$ denotes the last element of this sequence.}}}
\end{algorithm}

\begin{figure}
 \centering
 \includegraphics[width=2.0in]{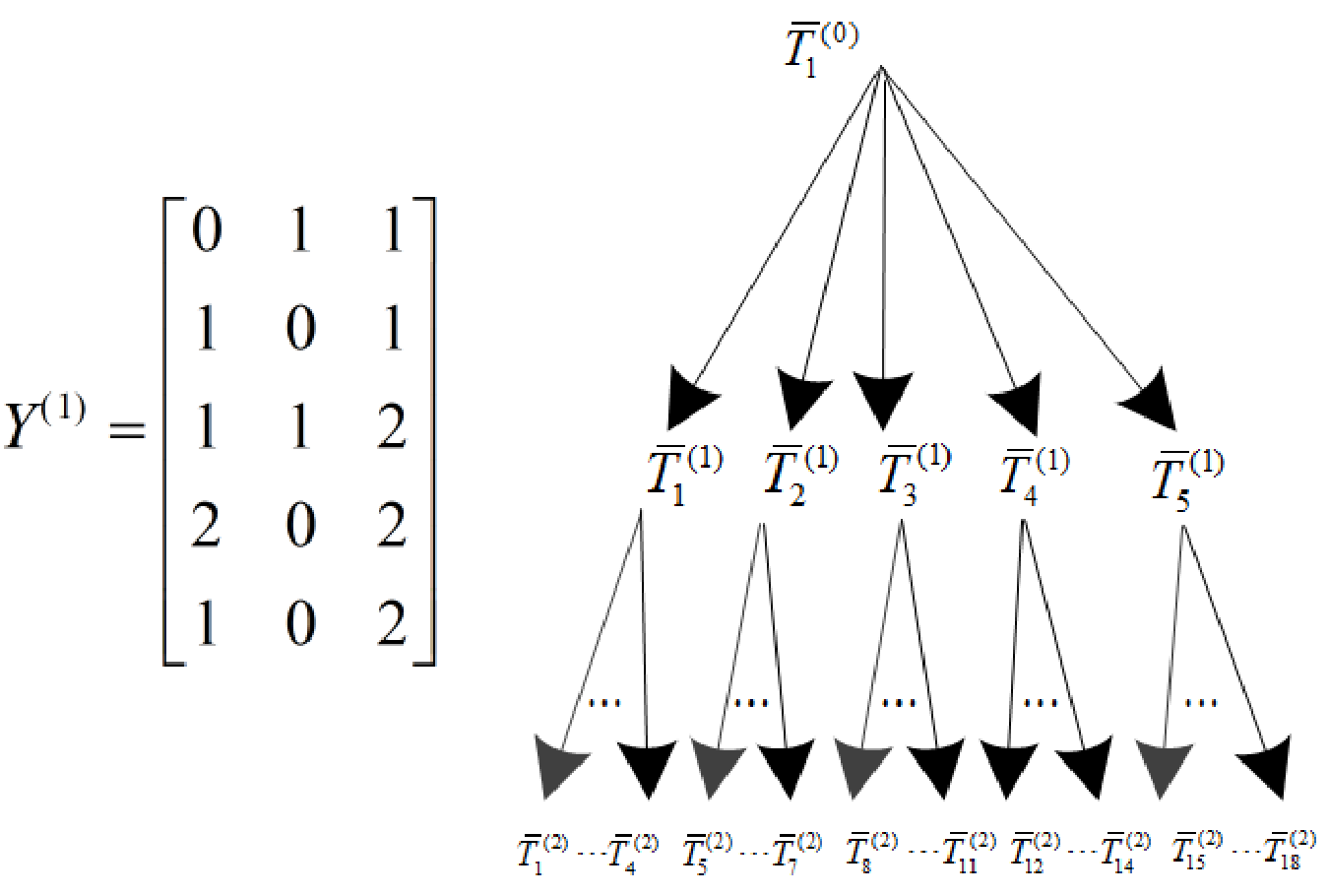}\\
 \caption{Illustration of the recursive tree in Algorithm \ref{alg1} associated with $Y^{(1)}$. In this example, $\bar T_{1}^{(0)}=\emptyset$, $\bar T^{(1)}_{i}=\{i\}$ for $i=1,3,5$, $\bar T^{(1)}_{2}=\bar T^{(1)}_{4}=\{2,4\}$, and $\bar T_{j}^{(2)}=\{1,2,3,4\}$ for $j=1,2,3,5,6,8,9,10,12,13$, $\bar T_{l}^{(2)}=\{2,4,5\}$ for $l=7,14,16,18$, $\bar T_{4}^{(2)}=\bar T_{15}^{(2)}=\{1,5\}$,  $\bar T_{11}^{(2)}=\bar T_{17}^{(2)}=\{3,5\}$.}\label{recursive tree}
 \end{figure}

In what follows, a polynomial time algorithm based on traversals over a recursive tree is presented. Here, `recursive' comes from the fact that this tree is {\emph{constructed recursively}}.  The pseudo code of this algorithm is given in Algorithm \ref{alg1}. The intuition behind Algorithm \ref{alg1} is that, for the $\kappa$th eigenvalue of $A$, $\kappa \in \{1,...,p\}$, instead of directly determining $r_\kappa$, we try to determine $r-r_{\kappa}$, which is the maximum number of rows of $Y^{(\kappa)}$ that fail to have full column rank. {However, since $r-r_{\kappa}$ could be as large as $r-1$, a naive combinatorial search needs to compute the ranks of at most ${\bf C}_r^{k_\kappa-1}+\cdots+{\bf C}_r^{r-1}$ submatrices, which again increases exponentially with $r$ even when $k_\kappa$ is bounded. The key idea to reduce the exponentially increasing complexity to the polynomially increasing one is leveraging the rank-one update rule expressed by (\ref{rank_one}).} Algorithm \ref{alg1} first builds a recursive tree and then searches the maximum return value among the leaves of this tree. An illustration of such a recursive tree is given in Fig. \ref{recursive tree}. To build the recursive tree for the $\kappa$th eigenvalue, the root is $\bar T^{(0)}_1=\emptyset$.
In the $\tau$th layer{\footnote{In this paper, the root of a tree is indexed as the $0$th layer.}}, $\tau=0,...,k_{\kappa}-1$, the set $T_i^{(\tau)}\subseteq V$ has the property that $Y^{(\kappa)}_{T_i^{(\tau)}}$ has rank $\tau$, $i=1,...,L_{\tau-1}$, where $L_{\tau-1}$ is the number of nodes in the $\tau$th layer. Then, $Y^{(\kappa)}_{\bar T_i^{(\tau)}}$ is obtained by adding the maximal rows of $Y^{(\kappa)}_{V\setminus T_i^{(\tau)}}$ to $Y^{(\kappa)}_{T_i^{(\tau)}}$ while its rank is preserved.  The collection of sets $\{\bar T_i^{(\tau)}|_{i=1}^{L_{\tau-1}}\}$ form the $\tau$th layer of the tree.  Next, each set ${\bar T_i^{(\tau)}}$ generates $|\Omega_i^{(\tau)}|$ sets, constituting $\{T_j^{(\tau+1)}|_{j=1}^{L_{\tau}}\}$ for the $(\tau+1)$th layer, where $\Omega_i^{(\tau)}$ is such set whose arbitrary element $j$ makes matrix $Y^{(\kappa)}_{\bar T_i^{(\tau)}\cup \{j\}}$ have rank $\tau +1$. After $\tau=k_{\kappa}-1$, the complete recursive tree is obtained, and one just needs to search the set with maximum cardinality among the leaves of this tree, i.e., $T_{\max }^{[\kappa ]}$, whose cardinality equals exactly $r-r_\kappa$.

\begin{theorem}\label{theorem 2}
Given a system $(A,C)$ with geometric multiplicities of eigenvalues of $A$ bounded by $\bar k$, $A\in \mathbb{R}^{n\times n}$, $C\in \mathbb{R}^{r\times n}$, Algorithm \ref{alg1} can determine the optimal solution to Problem 1 in time complexity at most $O({\bar k^2}{r^{\bar k+1}}n + r{n^3})$.
\end{theorem}

{\bf{Proof:}} As analyzed above, the main procedure of Algorithm~\ref{alg1} is to determine the maximum number of rows of $Y^{(\kappa)}$ that fail to have full column rank for each eigenvalue $\lambda_{\kappa}$, $\kappa=1,...,p$. To justify the recursive procedure, observe that for any $S\subseteq V$, $j\in V$,
\begin{equation}\label{rank_one}{\rm{rank}}Y_{S\bigcup {\{ j\} } }^{(\kappa )} = \left\{ \begin{array}{l}
{\rm{rank}}Y_S^{(\kappa )},{\rm if}{\kern 2pt}  \overline{{\rm span}}(Y_{\{j\}}^{(\kappa )}) \subseteq \overline{{\rm span}}(Y_S^{(\kappa )})\\
{\rm{rank}}Y_S^{(\kappa )} + 1,{\rm if}{\kern 2pt}   \overline{{\rm span}}(Y_{\{j\}}^{(\kappa )}) \not\subseteq \overline{{\rm span}}(Y_S^{(\kappa )})
\end{array} \right.\end{equation}
Hence, for the $\kappa$th eigenvalue, the recursive tree has depth exactly $k_{\kappa}$. Moreover, $\{\bar T_1^{({k_\kappa } - 1)}, \cdots ,\bar T_{L_{{k_\kappa}-2}}^{({k_\kappa } - 1)}\}$ index {\emph{all}} the submatrices formed by rows of $Y^{(\kappa)}$ with the property that: it has rank $k_{\kappa}-1$, and adding any rest rows from $Y^{(\kappa)}$ can make it have full column rank. The optimality of the solution returned from Algorithm \ref{alg1} then follows immediately.

In the recursive tree for the $\kappa$th eigenvalue, each parent node has at most $r$ child nodes. Thus, the $\tau$th layer has at most $r^{\tau}$ nodes, $\tau=0,...,k_{\kappa}-1$. Hence, the total nodes in that tree is at most $\sum\nolimits_{\tau = 0}^{{k_\kappa} - 1} {{r^\tau}}=(r^{k_\kappa}-1)/(r-1) \to O({r^{{k_\kappa}-1}})< {r^{{k_\kappa}}}$ ($r\ge 3$). For each node, the rank update procedure incurs $O(k_\kappa ^2r\cdot r)$ using the singular value decomposition (Lines 6-7 of Algorithm \ref{alg1}). To obtain $Y^{(\kappa)}|_{\kappa=1}^p$, it incurs computational complexity $O(rn^3)$. Note that $k_{\kappa}\le \bar k$ for $\kappa =1,...,p$, and $p\le n$. Hence, the total computational complexity of Algorithm \ref{alg1} is at most
 $O({\bar k^2}{r^{\bar k+1}}n + r{n^3})$. $\hfill\blacksquare$

{\begin{remark}
The above results show that the computational intractability of Problem 1 is essentially caused by the increase of eigenvalue geometric multiplicities of $A$, rather than that of dimensions of states or number of sensors. By contrast, from \cite{A.Ol2014Minimal,T2016On}, it is known that the MCPs discussed therein are NP-hard even when the involved systems have no repeated eigenvalues. To accelerate Algorithm \ref{alg1}, one could remove the repeated elements from $\{\bar T_i^{\tau}|_{i=1}^{L_{\tau-1}}\}$ after Line 6, and collect all the rows $J\subseteq \Omega_i^{(\tau)}$ that make $Y^{(\kappa)}_{\bar T_i^{(\tau)}\cup J}$ have rank $\tau+1$ as a single node for each $T_i^{(\tau)}$ in Lines 7-8, to {\emph{reduce the number of repeated nodes}} in the $\tau$th layer of the recursive tree.

%, which means the increase of system dimensions and number of actuators cause the exponential increase in the computation cost

%To accelerate Algorithm \ref{alg1}, instead of directly generating $|\Omega_i^{(\tau)}|$ nodes for each $T_i^{(\tau)}$ in Lines 7-8, one could collect all the rows $J\subseteq \Omega_i^{(\tau)}$ that make $Y^{(\kappa)}_{\bar T_i^{(\tau)}\cup J}$ have rank $\tau+1$ as a single node, to reduce the number of possibly repeated nodes in the $\tau$th layer of the recursive tree.
%This is perhaps the essential difference between the MCPs and the minSRO considered here. % (and the minARO problem). and remove the repeated elements in $\{\bar T_i^{\tau}|_{i=1}^{L_{\tua-1}}\}$ after Line 6,

\end{remark}
}

%{ Advantage of the proposed algorithm:} 1. dynamic rank-one update algorithm. the computational complexity is $O(mn(log\{m,n\})^2)$ if a column is added to $A\in {\mathbb R}^{m\times n}$. [Lemma 3.3 of \cite{Cheung-fastrank-2013}, \cite{Frandsen-dyrank-2009}]
%  { for more details, see the extended version of this paper (the updated version of [1]).}

\begin{remark} { It is remarkable that Algorithm \ref{alg1} essentially finds the optimal solutions by traversing a subset of feasible solutions $\{V\backslash \bar T_i^{(k_\kappa-1)}|_{i=1,...,L_{k_{\kappa}-2}}^{\kappa=1,...,p}\}$ to the original problem. Since it is an exact algorithm, without any restriction on the problem input, it cannot be guaranteed Algorithm \ref{alg1} will perform {\emph{significantly}} better than the exhaustive combinatorial search in the worst case. A {\emph{similar computational complexity}} to Algorithm \ref{alg1} can be obtained by an alternative two-stage procedure, which consists of, first enumerating all $k_{\kappa}-1$ independent rows of $Y^{(\kappa)}$ ($\kappa=1,...,p$), and then for each element, adding the maximal number of rows from $Y^{(\kappa)}$ to it without increasing its rank of $k_{\kappa}-1$. To ensure a polynomial time complexity with bounded $k_\kappa$, instead of enumerating over the $r-(k_\kappa-1)$ rows in the second stage, one should leverage the rank-one update property (\ref{rank_one}) and adopt a way similar to Line 6 of Algorithm \ref{alg1}. The only difference between those two procedures is that Algorithm \ref{alg1} uses a recursive (rank-one update) way to get a subset of feasible solutions. In contrast, the latter procedure directly skips to finding all rank-$(k_{\kappa}-1)$ bases. Note throughout Algorithm \ref{alg1}, the only computation involving the rank is to check whether a matrix can preserve its rank after adding a row (Line 6). Hence, some fast dynamic rank-one update algorithms \cite{Frandsen-dyrank-2009,Cheung-fastrank-2013} can be used to reduce the rank involved computations in each iteration. This may sometimes make Algorithm \ref{alg1} faster than the two-stage procedure. To be specific, directly implementing the two-stage procedure for Problem 1 incurs $O(\bar k^2r^{\bar k+1}n+rn^3)$ (the best bound is $O(\bar k^{w-1}r^{\bar k+1}n+rn^3)$), while Algorithm 1 using the dynamic rank-one update algorithm in \cite{Cheung-fastrank-2013} incurs $\tilde O(\bar kr^{\bar k+1}n+\bar k^w+rn^3)$, where $w=2.373$ is the matrix multiplication exponent.\footnote{$\tilde O(f(n))$ stands for $O(f(n)(log n)^k)$ for some positive constant $k$.} Moreover, the recursive rank-one update way may be more preferable in the scenario where partial sensors are forbidden to be removed, or where the eigenvalues $\lambda_i|_{i=1}^p$ are more easily accessible than the eigenbases $X_i|_{i=1}^p$ \cite{Geor1993Graph} (in which case one can alternatively minimize $|S|$ to make ${\bf col}\{A-\lambda_i I,C_{V\backslash S}\}$ column rank deficient, by Lemma \ref{lemma 1}). In these scenarios, most rows of $C$ could be collected to $\bar T^{({\tau})}_i$ in the early iterations of Algorithm \ref{alg1}, leading to a relatively small $\Omega_i^{({\tau})}$, thus a recursive tree with relatively small number of nodes.}
\end{remark}

\section{Computational Perspective of Problems 2 and 4}
Algorithm \ref{alg1} requires accurate eigenspace decomposition of the system matrices, which may encounter numerical instabilities or rounding errors for large-scale systems. Hence, studying the corresponding structured systems may be more practical.  In this section, we first prove that Problems 2 and 4 are NP-hard,
then provide the exact algorithm for Problem 2 under a corresponding assumption, and finally give a decomposition-based acceleration technique for this algorithm.
\vspace*{-3mm}
\subsection{Complexity of Problems 2 and 4}
\begin{theorem} \label{theorem2}
Problem 2 is NP-hard, { even when $n\!-\!r$, i.e., the difference between the numbers of states and sensors, is arbitrarily large.} % even when the number of states is arbitrarily larger than that of the sensors. %{ even when the number of states is arbitrarily larger than that of the sensors}. %
\end{theorem}

We first give some definitions and lemmas needed to proceed with our proof. Let ${\bf \bar 1}_{m\times n}$ be an $m\times n$ structured matrix whose entries are all free parameters. A set of columns of a structured matrix are said to be linearly dependent, if the generic rank of the matrix consisting of them is less than the number of their columns.  A circuit of an $m\times n$ matrix $M$, is a set $J\subseteq \{1,...,n\}$  such that ${\rm rank}(M_{[J]})=|J|-1$. The set of indices of any linearly dependent columns of $M$ contains at least one circuit of $M$. A clique of an undirected graph is a subgraph such that any two of its vertices are adjacent.  The $k$-clique problem is to determine whether there is a clique with $k$ vertices, which is NP-complete~\cite{Geor1993Graph}.

%\begin{lemma}
%Given a structured matrix $\bar M \in \{0,*\}^{n\times m}$ and an integer $r<n$, it is NP-complete to determine whether there exists a subset $J\subseteq \{1,...,m\}$ with $|J|=r$, such that ${\rm grank}(M_J)<r$.
%\end{lemma}
% {\bf{Proof:}}

\begin{lemma} \label{rank_less}
Let matrix $\bar M \in \{0,*\}^{n\times l}$. If ${\rm grank}([{\bf \bar 1}_{n\times m}, \bar M])<n$, where $m<n$, then ${\rm grank}(\bar M)<n-m$.
\end{lemma}

{\bf{Proof:}} Suppose ${\rm grank}(\bar M)=p\ge n-m$. Then, there exist two sets $J_1\subseteq \{1,...,n\}$ and $J_2\subseteq \{1,...,l\}$ with $|J_1|=|J_2|=p$, such that ${\rm grank}(\bar M_{J_1,J_2})\ge n-m$. Therefore, ${\rm grank}([{\bf \bar 1}_{n\times m}, \bar M])\ge  {\rm grank}(
\left[ {\begin{array}{*{20}c}
   {{\bf \bar 1}_{n - p,m} } & 0  \\
   {{\bf \bar 1}_{p,m} } & {\bar M_{J_1,J_2} }  \\
\end{array}} \right]
)=\min\{n-p,m\}+p=n-p+p= n$, causing a contradiction. $\hfill\blacksquare$

{\bf{Proof of Theorem \ref{theorem2}:}}
Our proof is a reduction from the NP-complete $k$-clique problem to one instance of Problem 2. We divide it into three steps.

Step 1). Let ${\cal G}=({\cal V}, {\cal E})$ be an undirected graph (without self-loop), and $|{\cal V}|=p$, $|{\cal E}|=r$. The incidence matrix of ${\cal G}$, denoted by  $In({\cal G})$, is an $r\times p$ matrix so that $[In({\cal G})]_{ij}=1$ if $v_j\in V(e_i)$ and $0$ otherwise, where $v_j\in {\cal V}$, $e_i\in \cal E$.
Construct the matrix $\tilde C\doteq \left[In({\cal G}),
{{{\bf \bar 1}_{r \times({\tiny{ \left( {\begin{array}{*{20}{c}}
k\\
2
\end{array}} \right)}} - k - 1)}}} \right]$, and let the output matrix $\bar C$ be $
\bar C=[{0_{r\times l}},\tilde C]$, where $k\ge 5$ (so that ${\tiny{ \left( {\begin{array}{*{20}{c}}
k\\
2
\end{array}} \right)}} - k - 1\ge 2$), `$1$' in $In(\cal G)$ represents a free entry, {and the integer $l\ge 0$ can be arbitrarily selected.} Let $q\doteq \tiny{ \left( {\begin{array}{*{20}{c}}
k\\
2
\end{array}} \right)}$. Without any losing of generality, assume that $p>k+1$. It demonstrates in  Page 53 of \cite{1985_spark} that: 1)  Every circuit $J$ of $\tilde C^{\intercal}$ must satisfy $|J|\ge q$. 2)  If $\cal G$ has a $k$-clique, then $\tilde C^{\intercal}$ has a circuit with cardinality $q$. 3) If $\tilde C^{\intercal}$ has a circuit with cardinality $q$, then $\cal G$ has an $k$-clique. %${\cal G}$ is connected and , which does not affect the NP-completeness of the $k$-clique problem associated with $\cal G$

From the aforementioned results, it can be declared that, $\tilde C^{\intercal}$ (as well as $\bar C^{\intercal}$) has $q$ linearly dependent columns, if and only if the digraph $\cal G$ has an $k$-clique. Indeed, from Fact 1), if $\tilde C^{\intercal}$ has $q$ linearly dependent columns,  then such columns must form a circuit of $\tilde C^{\intercal}$. Otherwise, $\tilde C^{\intercal}$ has a circuit with cardinality less than $q$. With Fact 3), the `only if' direction is obtained. The `if' direction comes directly from Fact 2).

Step 2). Construct the state transition matrix $\bar A$ as
\[
\bar A = \left[ {\begin{array}{*{20}c}
   {0_{q \times q} } & {0_{q\times (n-q)} }   \\
    {{\bf \bar 1}_{(n - q)\times q} } & {{\bf \bar 1}_{(n - q) \times (n - q)} }  \\
\end{array}} \right]
,\]where $n\doteq l+p+q-k-1>q$, {and $l$ could be selected so that $n-r$ can be arbitrarily large.} Construct the system digraph ${\cal D}(\bar A, \bar C)$, and let state vertices ${\cal X}=\{x_1,...,x_n\}$, output vertices ${\cal Y}=\{y_1,...,y_r\}$. From the construction, it is clear that every state vertex in ${\cal X}\backslash \{x_n\}$ is an in-neighbor of $x_n$, and $x_n$ is reachable to every output in $\cal Y$. Hence, all the $r$ outputs must be deleted so that at least one state vertex is output-unreachable. In other words, the minimum number of outputs (sensors) needed to be deleted to destroy the output-reachability of at least one state vertex equals $r$ (Lemma \ref{Lemma 1}).

Step 3). We declare that the optimal value of Problem 2 associated with $(\bar A,\bar C)$ is no more than $r-q$, if and only if there exist $q$ linearly dependent columns in $\bar C^{\intercal}$.  Denote the optimal solution to Problem 2 by $S_{opt}\subseteq V$.
For the one direction, suppose there are $q$ linearly dependent columns in $\bar C^{\intercal}$, and denote the set of its indices by $S_{ld}$. Let $S_{re}=V\backslash S_{ld}$.  It holds that $$\begin{array}{l}
{\rm{grank}}([{{\bar A}^{\intercal}},{{[\bar C^{\intercal}]}_{[V\backslash {S_{re}}]}}]) = {\rm{grank}}([{{\bar A}^{\intercal}},{{[\bar C^{\intercal}]}_{[{S_{ld}}]}}])\\
 \le {\rm{grank}}(\bar A) + {\rm{grank}}({{[\bar C^{\intercal}]}_{[{S_{ld}}]}}) < n - q + q = n.
\end{array}$$
As a result, $S_{re}$ is a feasible solution to Problem 2 associated with $(\bar A, \bar C)$. Since $|S_{re}|=r-q$, we have $|S_{opt}|\le r-q$.

For the other direction, suppose that $|S_{opt}|\le r-q$. Then, since the minimum number of sensors whose removal destroys the output-reachability condition equals $r$, it must hold that  ${\rm grank}([\bar A^{\intercal}, {[\bar C^{\intercal}]}_{[V\backslash S_{opt}]}])<n$. As ${\rm grank}(\bar A)=n-q$, from Lemma \ref{rank_less}, it is valid that
${\rm grank}({[\bar C^{\intercal}]}_{[V\backslash S_{opt}]})<q$. Noting that $|V\backslash S_{opt}|\ge q$, we have that the columns of $\bar C^{\intercal}$ indexed by $V\backslash S_{opt}$ are linearly dependent.

Since we have shown that verifying whether there are $q$ linearly dependent columns in $\bar C^{\intercal}$ (so in $\tilde C^{\intercal}$) is equivalent to verifying whether the digraph $\cal G$ has an $k$-clique, which is NP-complete, and the above-mentioned reduction is in polynomial time, it is concluded that Problem 2 is NP-hard.  $\hfill\blacksquare$

%While \cite{zhang2019minimal} has proven the NP-hardness of determining the minimal cost of inputs whose removal causes structural uncontrollability, it cannot be obtained by duality that the minSRO for structured systems is NP-hard, as in their proof different inputs may have different costs.  %Here to eliminate the ununiform cost used therein, we have adopted a different construction of $\bar A^{\intercal}$ and $\bar C^{\intercal}$.
The key challenge to prove Theorem \ref{theorem2} lies in constructing a system in which there is an explicit relationship in size between the number of sensors whose removal destroys the output-reachability condition and that whose removal destroys the matching condition.
{
In our previous work \cite{zhang2019minimal}, we have proved the NP-hardness of Problem 2 in the case with nonuniform sensor cost (by duality between controllability and observability), and indicated that it might be possible to extend the proof technique therein to show the computational complexity of Problem 2 with uniform cost (since the extension may not be so straightforward and much of the details is not given, the NP-hardness of Problem 2 has not been definitively declared in \cite{zhang2019minimal}). In that extension, by replicating one sensor $n-1$ times, where $n$ is the number of states of the constructed system instance, the number of sensors is $2n-1$. Here, we have provided a different technique yet a more {\emph{practical}} system instance to show the NP-hardness of Problem 2, in which any two sensors have {\emph{different}} structure configurations (i.e., they measure different combinations of state variables), and the number of state variables $n$ can be {\emph{arbitrarily larger}} than that of the sensors $r$.  Furthermore,  our proof technique is more general in the sense that it can be slightly modified to show the intractability of Problem 4.

% our proof technique can be slightly modified to show the NP-hardness of Problem 4. %Particularly, the difference $n-r$ can be arbitrarily large. %is more general in the sense that it can be ...
}

{
{\begin{theorem} Problem 4 is NP-hard. % Given a structured system $(\bar A,\bar C)$, it is NP-hard to determine the minimal number of state variables that need to be prevented from being directly accessible to the existing sensors (i.e., the number of nonzero columns of $\bar C$ that need to be zeroed out), such that the resulting system becomes structurally unobservable.
\end{theorem}}

{\bf{Proof:}} For space consideration, we give the sketch of our proof, which follows similar reasoning to the proof of Theorem \ref{theorem2} with slight modifications. Given an undirected graph ${\cal G}=({\cal V},{\cal E})$, let $In({\cal G})$ (the incidence matrix), $r\ (\doteq |{\cal E}|)$, $p\ (\doteq |{\cal V}|)$, $k\ (\ge 5)$, $q\ (\doteq {\tiny{ \left( {\begin{array}{*{20}{c}}
k\\
2
\end{array}} \right)}})$ and $\tilde C$ be defined in the same way as in the proof of Theorem \ref{theorem2}.
Construct a structured system $(\hat A, \hat C)$ as
$$\hat A\!=\!\left[\!\!
             \begin{array}{cc}
               0_{(q-1)\times (q-1)} & 0_{(q-1)\times (r+1-q)} \\
               {\bf \bar 1}_{(r+1-q)\times (q-1)} & {\bf \bar 1}_{(r+1-q)\times (r+1-q)} \\
             \end{array}
       \!  \!  \right]
, \hat C\!=\!\! \left[\!\!
                         \begin{array}{c}
                           [In({\cal G})]^{\intercal} \\
                           {\bf \bar 1}_{(q-k-2)\times r} \\
                         \end{array}
                       \!\!\right]
 .$$ Let $S_{opt}\subseteq \{1,...,r\}$ be the set of state variables with the minimum cardinality that need to be blocked from the existing sensors in system $(\hat A, \hat C)$, so that the resulting system becomes structurally uncontrollable. We will show that $|S_{opt}|\le r-q$, if and only if $\cal G$ has a $k$-clique. %Denote the optimal solution to this problem by $S_{opt}\subseteq \{1,...,r\}$.

For the one direction, suppose that $\cal G$ has a $k$-clique. Then, from Step 1) in the Proof of Theorem \ref{theorem2}, $\tilde C$ has $q$ linearly dependent rows, whose row index set is denoted by $S_{ld}$ ($|S_{ld}|=q$). Using Lemma \ref{rank_less} on $\tilde C_{S_{ld}}$, we get ${\rm grank} (\hat C_{[S_{ld}]})< q-1$, noting that ${\rm grank} (\tilde C_{S_{ld}})< q$ and $\tilde C_{S_{ld}}\equiv [[\hat C^{\intercal}]_{S_{ld}}, {\bf \bar 1}_{q\times 1}]$.
 Let $S_{re}\doteq V \backslash S_{ld}$, where $V=\{1,...,r\}$.
 Then, we have
${\rm grank}({\bf col}\{\hat A,\hat C^{S_{ld}}\})\le {\rm grank}(\hat A)+ {\rm grank}({\hat C}_{[S_{ld}]})< r+1-q+q-1=r,$
which means $(\hat A,\hat C^{V\backslash S_{re}})$ is structurally uncontrollable by Lemma \ref{lemma 1}. Since $|S_{re}|=r-q$, we get $|S_{opt}|\le |S_{re}|=r-q$.

%$$\begin{array}{c}{\rm grank}({\bf col}\{\hat A,\hat C^{\backslash S_{re}}\})\le {\rm grank}(\hat A)+ {\rm grank}({\hat C}_{[S_{ld}]})\\< r+1-q+q-1=r,\end{array}$$ %%Recall that $\hat C^{S_{re}}$ is the matrix obtained from $\hat C$ by zeroing out its columns indexed by $S_{re}$.

For the other direction, suppose $|S_{opt}|\le r-q$. Let ${\cal X}={\cal X}_1\cup {\cal X}_2$ be the state vertex set of the digraph ${\cal D}(\hat A,\hat C)$, with ${\cal X}_1=\{x_1,...,x_{q-1}\}$, ${\cal X}_2=\{x_{q},...,x_{r}\}$. Then, it is easy to see that, any two vertices in ${\cal X}_2$ are reachable from each other, and every vertex in ${\cal X}_1$ is reachable to arbitrary vertex in ${\cal X}_2$. Therefore, all vertices in ${\cal X}_2$ must be blocked from the corresponding sensors, so that the out-reachability condition is violated. As $|{\cal X}_2|=r-q+1> r-q$ and $|S_{opt}|\le r-q$, to make $(\hat A, \hat C^{V \backslash S_{opt}})$ structurally uncontrollable, it must hold that
${\rm grank}({\bf col}\{\hat A,\hat C^{V \backslash S_{opt}}\})={\rm grank}([\hat A^{\intercal},[\hat C^{V \backslash S_{opt}}]^{\intercal}]) < r$. Again using
Lemma \ref{rank_less} on $[\hat A^{\intercal},[\hat C^{V \backslash S_{opt}}]^{\intercal}]$, we get ${\rm grank}(\hat C_{[{V \backslash S_{opt}}]})={\rm grank}(\hat C^{V \backslash S_{opt}})< r-(r+1-q)=q-1$. Therefore, ${\rm grank}([[\hat C^{\intercal}]_{V\backslash S_{opt}}, {\bf \bar 1}_{|V\backslash S_{opt} |\times 1}])\le {\rm grank}([\hat C^{\intercal}]_{V\backslash S_{opt}})+1<q-1+1=q$. Noting that $[[\hat C^{\intercal}]_{V\backslash S_{opt}}, {\bf \bar 1}_{|V\backslash S_{opt} |\times 1}]$ is exactly $\tilde C_{V\backslash S_{opt}}$ and $|V \backslash S_{opt}|\ge q$, we therefore obtain $\tilde C^{\intercal}$ has $q$ linearly dependent columns indexed by $V\backslash S_{opt}$. Immediately, from Step 1) in the proof of Theorem \ref{theorem2}, $\cal G$ has a $k$-clique.

Since the $k$-clique problem is NP-complete, we conclude that determining $|S_{opt}|$ is NP-hard. $\hfill\blacksquare$
}

\subsection{The CCSM Structure of Problem 2}
Similar to Section \ref{CCSM-sec-pro1}, Problem 2 also has a CCSM structure. To see this, we first resort to the strongly connected component (SCC) decomposition to determine the minimal sensors (denoted by $J^y_{re}$) whose removal destroys the output-reachability condition.

 A digraph is said to be strongly connected if there is a path from every one of its vertices to the other one. An SCC of $\cal G$ is a subgraph that is strongly connected and no additional edges or vertices from $\cal G$ can be included in this subgraph without breaking its property of being strongly connected. We call an SCC the sink SCC if there is no outgoing edge from this SCC to any other SCCs.
Suppose there are $l$ sink SCCs in ${\cal D}(\bar A)=({\cal X},{\cal E}_{{\cal X},{\cal X}})$. For each sink SCC ${\cal N}_i$, denote by $N({\cal N}_i)$ the out-neighbors of ${\cal N}_i$ in ${\cal D}(\bar A, \bar C)$ (hence $N({\cal N}_i)\subseteq {\cal Y}$). Define $J^y_{i}=\{j: y_j\in N({\cal N}_i)\}$; that is, $J^y_{i}$ is the set of sensors that are out-neighbors of ${\cal N}_i$. Then, it is easy to verify that, the minimum number $|J^y_{re}|$ of sensors  whose removal destroys the output-reachability condition equals $\min\nolimits_{1\le i\le l} |J^y_i|$.

According to Lemma \ref{Lemma 1} and similar to (\ref{ccsm1}), Problem 2 can be solved by invoking the following problem for at most $O(log_2|J^y_{re}|)$ times until the objective value is less than $n$.
\begin{equation}\label{ccsm-2}
\min\limits_{S \subseteq V, |S|\le |J^y_{re}|}  {\rm grank}\ {\bf col }\{\bar A, \bar C_{V\backslash S}\}
\end{equation}Again, (\ref{ccsm-2}) belongs to the CCSM problem, as the generic rank of sub-rows of a matrix is submodular on its row set (see \citep[Sec 2.3.3]{Murota_Book}). This indicates that, Problem 2 might be hard to approximate.

%The above problem belongs to the cardinality-constraint submodular function minimization (SFM).
\subsection{Algorithm with Bounded Matching Deficiencies}\label{alg-sec-stru}
We will call the number of vertices of ${\mathcal R}(\bar A)$ that are not matched in any maximum matching of ${\cal B}(\bar A)=({\cal R}(\bar A), {\cal C}(\bar A), {\cal E}(\bar A))$, as the {\emph{matching deficiency}}, which is equal to $n-{\rm grank}(\bar A)$.

\begin{assumption}\label{asp1}
The matching deficiency of $\bar A$ is bounded by some fixed constant $\bar k_s$ as $n$ and $r$ increase.
\end{assumption}

\begin{remark}
Assumption \ref{asp1} is satisfied by many practical systems. For example, in consensus networks or multi-agent systems \cite{multi_agent_2012}, every state vertex often has a self-loop, under which circumstance ${\mathcal B}(\bar A)$  has a zero matching deficiency. Moreover, connected ER random networks can be controlled using only one dedicated input with probability $1$ \cite{o2016conjecture}, which means the associated matching deficiencies are no more than $1$ \cite{Y.Y.2011Controllability} with probability~$1$. { The scenarios where the
matching deficiency is bounded by a constant which is neither zero nor one may arise in weakly-connected directed graphs that have a bounded (by a known constant) number of roots or leaves, and graphs with special layered structures, such as being the Cartesian products of factor graphs \cite{hammack2011handbook}.  }  %(and for most of time exactly zero) \cite{Y.Y.2011Controllability}
\end{remark}

%\subsubsection{Algorithm for Problem 2}
To determine the minimal sensors whose removal destroys the matching condition (denoted by $J^y_{ma}$), the similar idea to Algorithm \ref{alg1}, by constructing a recursive tree, is utilized. This tree has a similar structure to Fig. \ref{recursive tree}; see the pseudo code in Algorithm \ref{alg3}.  In the $\tau$th layer, each node $\bar T_i^{(\tau)}$ indexes a subset of sensors such that ${\rm grank}({\bf col}\{\bar A, \bar C_{\bar T_i^{(\tau)}}\})={\rm grank}(\bar A)+\tau$, and no other element from $V\backslash \bar T_i^{(\tau)}$ can be added to $\bar T_i^{(\tau)}$ without increasing ${\rm grank}({\bf col}\{\bar A, \bar C_{\bar T_i^{(\tau)}}\})$. Hence, in the $(n-{\rm grank}(\bar A)-1)$th layer, the leaf with the maximum cardinality is the complementary set to the minimal rows of $\bar C$ whose removal from ${\bf col}\{\bar A, \bar C\}$ makes the resulting matrix have generic rank $n-1$. This means the total layers of the recursive tree is exactly $n-{\rm grank}(\bar A)$.

%Hence, in the $\tau_m\doteq (n-{\rm grank}(\bar A)-1)$th layer, each node $\bar T_i^{\tau_m}$ corresponds to a subset of sensors such that ${\rm grank}({\bf col}\{\bar A, \bar C_{\bar T_{\tau_m}^{\tau}}\})=n-1$, and adding any element from $V\backslash \bar T_i^{\tau_m}$ to $\bar T_i^{\tau_m}$ will make the associated generic rank reaching $n$. Hence, the total layers of the recursive tree is $n-{\rm grank}(\bar A)$.

% Then, the $1$st layer consists of $r$ nodes. In the $i$th layer, there are at most $r^{i}$ nodes, $i=0,...,n-{\rm grank}(\bar A)-1$. In the ($n-{\rm grank}(\bar A)-1$)th layer, the leaf with the maximum cardinality is the set corresponding to the maximum rows of $\bar C$ whose removal from ${\bf col}\{\bar A, \bar C\}$ makes the resulting matrix having generic rank $n-1$.

\begin{theorem} \label{theo5}
Given a system $(\bar A, \bar C)$, where $\bar A \in \{0,*\}^{n\times n}$, $\bar C \in \{0,*\}^{r \times n}$, under Assumption \ref{asp1}, Algorithm \ref{alg3} provides an optimal solution to Problem 2 with time complexity at most $O((r+n)^{2.373}r^{\bar k_s})$.
\end{theorem}

{\bf{Proof}}: Based on the above arguments, to prove the optimality, we only need to {prove $\Delta$ defined in Line 6 of Algorithm \ref{alg3} is the maximal subset of rows of $\bar C_{V\backslash T_i^\tau}$ whose addition to ${\bf col}\{\bar A, \bar C_{T_i^\tau}\}$ can preserve its generic rank. Upon letting $\bar M_S\doteq {\bf col}\{\bar A, \bar C_S\}$, this can be justified from the following observation:  for any $S\subseteq V$, $j_1,j_2\in V$, if ${\rm grank} \bar M_{S}= {\rm grank}\bar M_{S\cup \{j_1\}}={\rm grank} \bar M_{S\cup \{j_2\}}$, then
${\rm grank} (\bar M_{S\cup \{j_1,j_2\}})\le {\rm grank} (\bar M_{S\cup \{j_1\}})+{\rm grank} (\bar M_{S\cup \{j_2\}})- {\rm grank} (\bar M_{S})= {\rm grank} (\bar M_S)$, leading to ${\rm grank} (\bar M_{S\cup \{j_1,j_2\}})= {\rm grank} (\bar M_{S})$, where the first inequality is due to the submodularity of the generic rank function of a matrix w.r.t. its rows \citep[Sec 2.3.3]{Murota_Book}.} This is similar to the rank-one update property (\ref{rank_one}).

% The proof for optimality of Algorithm \ref{alg3} follows the above arguments, and is similar to that of Theorem \ref{theorem 2}.
As for complexity, the SCC decomposition along with the step to obtain $J^y_{re}$ has complexity at most $O(n^2+nr)$ \cite{Geor1993Graph}. In the $i$th layer, there is at most $r^{i}$ nodes, $i=0,...,n-{\rm grank}(\bar A)-1$. Hence, the number of total nodes in the recursive tree is at most $\sum \nolimits_{\tau=0}^{\bar k_s-1} r^{\tau}= (r^{\bar k_s}-1)/(r-1)\to O(r^{\bar k_s-1})$ ($r\ge 3$). In each node, the operation updating $\bar T_i^{\tau}$ can be done by calling the maximum matching algorithm on bipartite graphs for at most $r$ times, which incurs complexity $O((r+n)^{2.373}r)$ in the worst case \cite{Geor1993Graph}. Hence, the total time complexity of Algorithm \ref{alg3} is $O(n^2+rn+(r+n)^{2.373}r^{\bar k_s})$, i.e, $O((r+n)^{2.373}r^{\bar k_s})$. $\hfill\blacksquare$

\begin{algorithm} % 算法开始
 {\small {{{
\caption{: Algorithm for Problem 2} %算法的题目
\label{alg3} %算法的标签
\begin{algorithmic}[1] %此处的[1]控制一下算法中的每句前面都有标号
\REQUIRE Parameters $(\bar A,\bar C)$ of a structurally controllable system  % a permutation $g=\{k_1,...,k_p\}$ of the integer set $\{1,...,p\}$
  %输入条件(此处的REQUIRE默认关键字为Require，在上面已自定义为Input)
\ENSURE   The optimal solution $J_{opt}$ to Problem 2  %输出结果(此处的ENSURE 默认关键字为Ensure在上面已自定义为Output)
\STATE  Find all the sink SCCs $\{{\cal N}_i|_{i=1}^l\}$ of ${\cal D}(\bar A)$. Let $J^y_{i}=\{j: y_j\in N({\cal N}_i)\}$. $J^y_{re}=\arg \min\nolimits_{1\le i\le l} |J^y_i|$.
\STATE  Initialize $\tau {\rm{ = }}0$, $T_1^{0} = \emptyset$, $L_{-1}=1$.
 \IF {{$n-{\rm grank}(\bar A)>0$}}
\WHILE{$\tau  < n-{\rm grank}(\bar A)$}
%\COMMENT{// We use $\{c^{\tau}_0,c^{\tau}_1,...,c^{\tau}_{\rm end}\}$ to denote a sequence of sorted integers with updated lengths, where $c^{\tau}_{\rm end}$ denotes the last element of this sequence.}
\FOR{$i=1$ to $L_{\tau-1}$}
%\COMMENT{We use $\{c^{\tau}_0,c^{\tau}_1,...,c^{\tau}_{\rm end}\}$ to denote a sequence of sorted integers with updated lengths, where $c^{\tau}_{\rm end}$ denotes the last element of this sequence.}
\STATE $\bar T_i^{\tau}=T_i^{\tau}\cup \Delta $, with $\Delta \!\doteq \! \{j: {\rm{grank}}({\bf{col}}\{\bar A,{{\bar C}_{T_i^{\tau} \cup \{ j\} }}\} )={\rm{grank}}({\bf{col}}\{ \bar A,{{\bar C}_{T_i^{\tau}}}\} ),j \in V\backslash T_i^{\tau}\}.$
%\STATE $\begin{array}{l}
%\bar T_i^{\tau} = T_i^{\tau} \cup \{ j:{\rm{grank}}({\bf{col}}\{ \bar A,{{\bar C}_{T_i^{\tau} \cup \{ j\} }}\} )\\
%{\kern 1pt} {\kern 1pt} {\kern 1pt} {\kern 1pt} {\kern 1pt} {\kern 1pt} {\kern 1pt} {\kern 1pt} {\kern 1pt} {\kern 1pt} {\kern 1pt} {\kern 1pt} {\kern 1pt} {\kern 1pt} {\kern 1pt} {\kern 1pt} {\kern 1pt} {\kern 1pt} {\kern 1pt} {\kern 1pt} {\kern 1pt} {\kern 1pt} {\kern 1pt} {\kern 1pt} {\kern 1pt} {\kern 1pt} {\kern 1pt} {\kern 1pt} {\kern 1pt} {\kern 1pt} {\kern 1pt} {\kern 1pt} {\kern 1pt} {\kern 1pt} {\kern 1pt} {\kern 1pt} {\kern 1pt} {\kern 1pt} {\kern 1pt} {\kern 1pt} {\kern 1pt} {\kern 1pt} {\kern 1pt} {\kern 1pt} {\kern 1pt} {\kern 1pt} {\kern 1pt} {\kern 1pt} {\kern 1pt} {\kern 1pt} {\kern 1pt} {\kern 1pt} {\kern 1pt} {\kern 1pt} {\kern 1pt} {\kern 1pt} {\kern 1pt} {\kern 1pt} {\kern 1pt} {\kern 1pt} {\kern 1pt} {\kern 1pt} {\kern 1pt}  = {\rm{grank}}({\bf{col}}\{ \bar A,{{\bar C}_{T_i^{\tau}}}\} ),j \in V\backslash T_i^{\tau}\}.
%\end{array}$
 \STATE    $\Omega _i^{\tau} =V\backslash \bar T_i^{\tau}$. Let $c_i^{\tau} = \sum\nolimits_{l = 1}^i {{\rm{|}}\Omega _l^{\tau}{\rm{|}}}$, $c_0^{\tau}{\rm{ = }}0$.
 \STATE Rewrite $\Omega _i^{\tau} \doteq \{ {j_1},...,{j_{|\Omega _i^{\tau}|}}\}$ and let $T_{c_{i - 1}^{\tau} + q}^{\tau  + 1}{\rm{ = }}\bar T_i^{\tau}\bigcup {j_q}$ for $q = 1,...,|\Omega _i^{\tau}|$.
%find the minimal {\rm span}ning forest of digraph $({\mathcal{Y}}\bigcup {{\mathcal{X}}, E(\mathcal{T})\cup \mathcal{M})$ with edge cost $C''$ rooted in $\mathcal{Y}$, given by $\mathcal{T}'$;
\ENDFOR
\STATE $L_{\tau}=c^{\tau}_{L_{\tau-1}}$; $\tau = \tau +1$.
%\COMMENT{THIS IS}
\ENDWHILE
 \STATE $T_{\max }= \arg \mathop {\max } \left\{ |\bar T_1^{{\tau_m } - 1}|, \cdots ,|\bar T^{\tau_m-1}_{L_{\tau_m-2}}|\right\}$, ${J^y_{ma}} = V\backslash T_{\max }$, with $\tau_m\doteq n-{\rm grank}(\bar A)$.
\ELSE
\STATE  {{ $J^y_{ma}=V$.}}
\ENDIF
 \STATE Return $J_{opt}=\arg \min \{|J^y_{ma}|, |J^y_{re}|\}$.
\end{algorithmic}}}
}}
\end{algorithm}

{\begin{remark}It can be seen both Line 6 of Algorithms \ref{alg1} and \ref{alg3}, which result from the rank-one update property of the involved rank function, are crucial to the polynomial time complexity  under the corresponding assumptions. However, the similar property does not hold for the general minSBOs. Hence, even under the same constant bound assumptions, it seems unclear whether they can be solved in polynomial time.
\end{remark}}

%\subsection{Systems with multiple eigenvalues}
%A more complicated case is when the system dynamics $A$ has multiple eigenvalues. Notice that, for a system with eigenvalues whose geometric multiplicity larger than one, the hitting set condition is no matter sufficient. We first introduce the necessary and sufficient condition for observability of general systems, which involves both graphic-algebraic elements [ACC paper].
%
%(a graphic-algebraic condition for observability of general STMs)
%
%Using the above criterion, in the following, we show that it is generally NP-hard to determine the minimum number of sensors to destroy observability even in the distributed control case.
%
%The reduction is from the matching preclusion number

\subsection{Acceleration of Algorithm \ref{alg3}}
The main procedure dominating the complexity of Algorithm~\ref{alg3} is to obtain the minimal sensors $J^y_{ma}$ whose removal destroys the matching condition. A method to accelerate this procedure is first to do Dulmage-Mendelsohn decomposition (DM-decomposition) on $\bar A$. Then implement Lines 3-15 of Algorithm \ref{alg3} on the reduced system, which usually has a dimension much less than that of the original system.

\begin{definition}[DM-decomposition,\cite{Murota_Book}]
Let $\bar M\in \{0,*\}^{m\times n}$. There exist two permutation matrices $P_r\in {\mathbb R}^{m\times m}$ and $P_c\in {\mathbb R}^{n\times n}$, such that
{\small\[P_r\bar MP_c=\left[
  \begin{array}{cccc}
    \bar M_{11} & \bar M_{12} & \bar M_{13} & \bar M_{14}\\
    0 & 0 & \bar M_{23} & \bar M_{24} \\
    0 & 0 & 0 & \bar M_{34} \\
    0 & 0 & 0 & \bar M_{44} \\
  \end{array}
\right],\]}where $\bar M_{12}$, $\bar M_{23}$, and $\bar M_{34}$ are square with zero-free diagonals. The columns of $\bar M_{11}$ (or the rows of $\bar M_{44}$) correspond to the unmatched vertices of at least one maximum matching in ${\cal B}(\bar M)$, and the rows (or columns) of $\bar M_{23}$ correspond to the vertices that appear in every maximum matching in ${\cal B}(\bar M)$. The generic rank of $\bar M$ equals $m$ minus the number of rows of $\bar M_{44}$.
\end{definition}

Using the DM-decomposition, suppose there are two permutation matrices $P_r$ and $P_c$ with compatible dimensions, such that
{\small\[P_r [\bar  A^{\intercal},\bar C^{\intercal}]\left[ {\begin{array}{*{20}c}
   {P_c } & {}  \\
   {} & I  \\
\end{array}} \right] = \left[ {\begin{array}{*{20}c}
   {\bar  A_{11} } & {\bar  A_{12} } & {\bar  A_{13} } & {\bar  A_{14} } &\vline &  {\bar B_1 }  \\
   {0} & {0} & {\bar  A_{23} } & {\bar  A_{24} } &\vline &  {\bar B_2 }  \\
   {0} & {0} & {0} & {\bar  A_{34} } &\vline &  {\bar B_3 }  \\
   {0} & {0} & {0} & {\bar  A_{44} } &\vline &  {\bar B_4 }  \\
\end{array}} \right],\]}where the left block of the right-hand side matrix of the above-mentioned formula is the DM-decomposition of $\bar A^{\intercal}$.
Since {\tiny{$\left[
        \begin{array}{cc}
          \bar A_{12} & \bar A_{13} \\
          0 & \bar A_{23} \\
        \end{array}
      \right]$}} has full generic rank, it can be validated that
\[{\rm grank}([\bar A^{\intercal},\bar C^{\intercal}]) = \left| {{\mathcal C}(\left[ {\begin{array}{*{20}c}
   {\bar A_{12} } & {\bar A_{13} }  \\
   0 & {\bar A_{23} }  \\
\end{array}} \right])} \right| + {\rm grank}(\left[ {\begin{array}{*{20}c}
   {\bar A_{34} } & {\bar B_3 }  \\
   {\bar A_{44} } & {\bar B_4 }  \\
\end{array}} \right])
\] %in Algorithm \ref{alg3}) whose removal destroy the matching condition
where ${\mathcal C}(\bullet)$ denotes the set of column indices of a matrix. Hence, to determine $J^y_{ma}$, we just need to identify the minimum number of columns of ${\bf {col}}\{\bar B_3,\bar B_4\}$ whose deletion makes $\left[{\tiny{ {\begin{array}{*{20}c}
   {\bar A_{34} } & {\bar B_3 }  \\
   {\bar A_{44} } & {\bar B_4 }  \\
\end{array}}}} \right]$ row generic rank deficient. As ${\bf {col}}\{\bar B_3,\bar B_4\}$ is obtained by permutating columns of $\bar C$, the set of indices of the aforementioned columns exactly correspond to $J^y_{ma}$. Hence, the reduced system $$\left({ \left[
                                                                                               \begin{array}{cc}
                                                                                                \bar A^{\intercal}_{34} & \bar A^{\intercal}_{44} \\
                                                                                                 0 & 0 \\
                                                                                               \end{array}
                                                                                             \right]
},[\bar B^{\intercal}_3,\bar B^{\intercal}_4] \right)$$ has the same solution as that of the original system $(\bar A,\bar C)$ on $J^y_{ma}$, where zero matrix is added to $[\bar A^{\intercal}_{34}, \bar A^{\intercal}_{44}]$ to square it. For a sparse matrix $\bar A$,  $|{\mathcal C}([\bar A^{\intercal}_{34}, \bar A^{\intercal}_{44}])|$ is usually much less than $|{\mathcal C}(\bar A)|$. Hence, it is perhaps safe to expect that the DM-decomposition can accelerate Algorithm \ref{alg3}.

%\section{observability robustness for structured systems}
%The observability robustness for a structured system with distributed control is NP-hard when actuator costs are polynomially bounded by the system size.
%
%The observability robustness for a structured system under distributed control with unit actuator costs is NP-hard.
\section{Numerical Results}\label{numer_example}
\begin{figure}
\centering
\subfigure[$c=3$]
{\includegraphics[width=1.75in, angle=0]{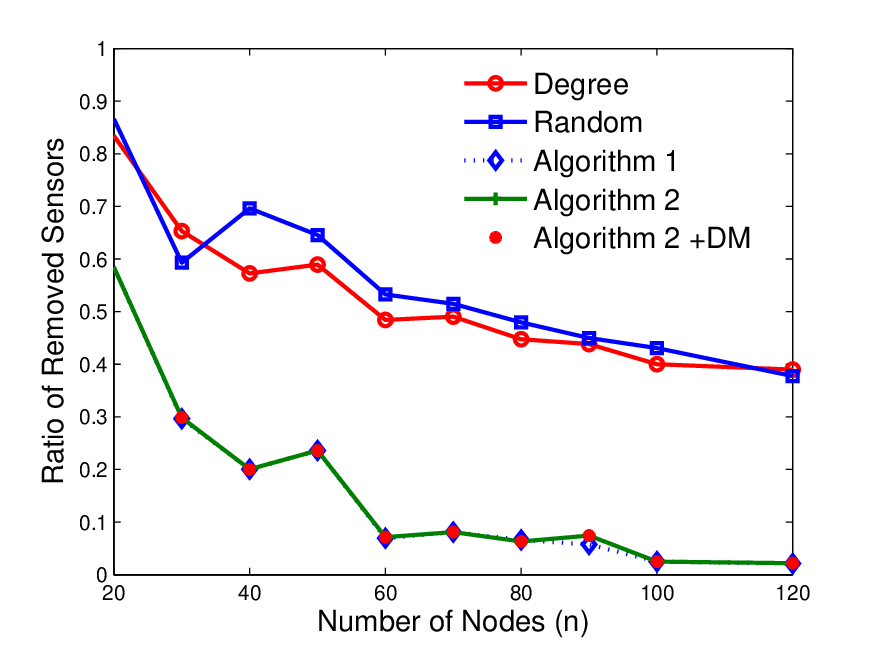}}
\subfigure[$c=4$]
{\includegraphics[width=1.75in, angle=0]{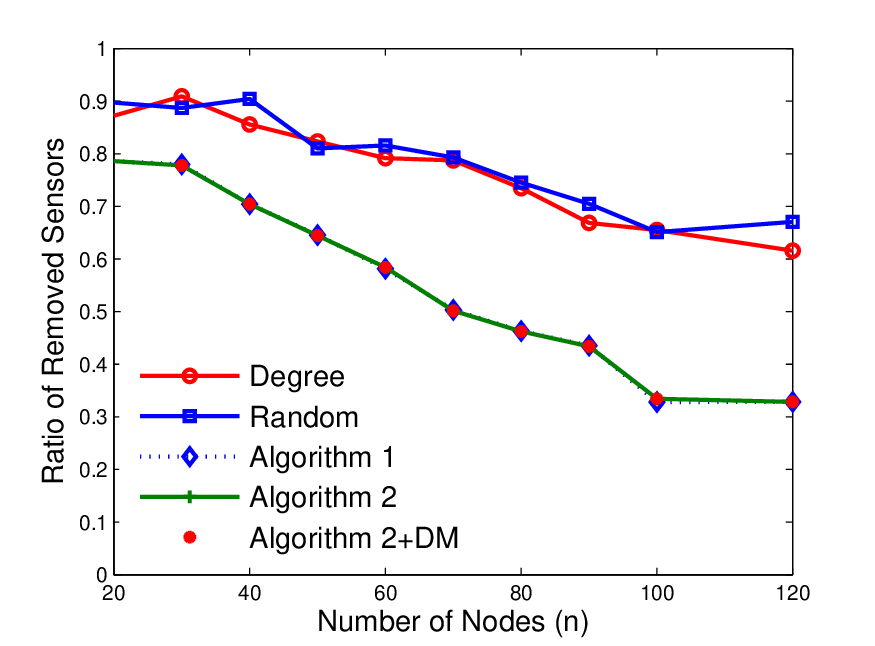}}
\caption{Experiments on ER random networks. Note that the vertical axis is the ratio between the number of removed sensors and that of the total sensors.}\label{figureSensorRemoval}
\end{figure}

We implement some numerical experiments to validate the effectiveness of our proposed algorithms. In our simulations, we generate ER random networks with $n$ nodes and use their adjacent matrices as the state transition matrices, $n$ ranging from $20$ to $120$. For each $n$, the probability of the existence of a directed edge between two nodes is set to be $\frac{c}{n}$, which means that the average degree (the sum of in-degree and out-degree) is $2c$ \cite{Geor1993Graph}. The weight of each directed edge is randomly generated in the interval $[-1,1]$. For each $n$, randomly generate $40$ observable systems (to save experiment time, all these systems have maximum geometric multiplicities less than $5$), where $40\%$ of the nodes are randomly selected to be measured by dedicated sensors (i.e., each sensor measures exactly one node). Five methods are adopted to identify a set of sensors whose removal destroys (structural) observability, the first of which is selecting the measured nodes in the descending order of their in-degrees, the second of which is randomly selecting the measured nodes, and the rest three of which are respectively based on Algorithm \ref{alg1}, Algorithm \ref{alg3} and Algorithm \ref{alg3} accelerated by the
DM-decomposition. The results are shown in Fig. \ref{figureSensorRemoval}.

%\begin{figure}
%\centering
%\begin{subfigure}{1.5\textwidth}
%  \centering
%  \includegraphics[width=1.4\linewidth]{SensorRemoval_NumberNodes_AvgDegree3Input4_30.eps}
%  \caption{A subfigure}
%  \label{fig:sub1}
%\end{subfigure}%
%\begin{subfigure}{1.5\textwidth}
%  \centering
%  \includegraphics[width=1.4\linewidth]{SensorRemoval_NumberNodes_AvgDegree4Input4_30.eps}
%  \caption{A subfigure}
%  \label{fig:sub2}
%\end{subfigure}
%\caption{A figure with two subfigures}
%\label{fig:test}
%\end{figure}
%\begin{figure}
%\centering
%\parbox{5cm}{
%\includegraphics[width=5cm]{SensorRemoval_NumberNodes_AvgDegree3Input4_30.eps}
%\caption{First.}
%\label{fig:2figsA}}
%\qquad
%\begin{minipage}{5cm}
%\includegraphics[width=5cm]{SensorRemoval_NumberNodes_AvgDegree4Input4_30.eps}
%\caption{Second.}
%\label{fig:2figsB}
%\end{minipage}
%\end{figure}

%\begin{algorithm} % 算法开始
%{\small  {{{{
%\caption{: Algorithm for Problem 4} %算法的题目
%\label{alg4} %算法的标签
%\begin{algorithmic}[1] %此处的[1]控制一下算法中的每句前面都有标号
%\REQUIRE System parameters $(\bar A,\bar C)$ %, a permutation $g=\{k_1,...,k_p\}$ of the integer set $\{1,...,p\}$
%  %输入条件(此处的REQUIRE默认关键字为Require，在上面已自定义为Input)
%\ENSURE   The optimal solution $\hat J_{opt}$ to Problem 4  %输出结果(此处的ENSURE 默认关键字为Ensure在上面已自定义为Output)
%\STATE Determine $\hat J^y_{re}$ using SCC-decomposition.
%\STATE Determine $\hat J^y_{ma}$ by running the procedure in {\bf{Lines 2-15 of Algorithm \ref{alg3}}}, in which one should replace Line 6 with $\bar T_i^{\tau}=T_i^{\tau}\cup \left\{j: {\rm grank}({\bf col}\{\bar A, \bar C^{T_i^{(\tau)}\cup \{j\}}\})= {\rm grank}({\bf col}\{\bar A, \bar C^{T_i^{(\tau)}}\}), j\in \hat V\backslash T_i^{\tau} \right\}$, and $V\to \hat V$.
% \STATE Return $\hat J_{opt}=\arg \min \{|\hat J^y_{ma}|, |\hat J^y_{re}|\}$.
%\end{algorithmic}}}
%}}}
%\end{algorithm}

Fig. \ref{figureSensorRemoval} shows that for each fixed $n$, the ratio of removed sensors returned by Algorithm \ref{alg1} is much less than that of the degree based and the random selection methods. This is reasonable, as Algorithm \ref{alg1} returns the optimal solutions. Moreover,  Algorithms \ref{alg1}, \ref{alg3} and the DM-decomposition accelerated Algorithm \ref{alg3} return almost the same values for each fixed $n$.\footnote{A deep insight shows that for some $n$, Algorithm \ref{alg1} returns values slightly smaller than those of Algorithm \ref{alg3}. This is consistent with the fact that for a given $(A, C)$, the optimal value of Problem 2 is a theoretical upper bound of that of Problem 1.} This indicates that, when the exact values of system matrices are not accessible, the structured system based Algorithm \ref{alg3} could be an alternative for the observability robustness evaluation. Comparing Fig. \ref{figureSensorRemoval}(a) and (b), it seems that networks with stronger connectivity (i.e., bigger average degrees) tend to have better observability robustness under sensor attacks.

\section{Conclusions}
%%Motivated by the question whether it is possible to uniquely reconstruct states of an LTI system under cardinality-constrained sensor attacks,
%The problem of determining the minimal number of sensors whose removal destroys system observability is considered. The dual of this problem is also the inverse to the minimal controllability problems \cite{A.Ol2014Minimal}. It is shown that this problem is NP-hard, both for a numerically specific system, and for a structured system. This confirms the conjecture in secure state estimation under adversarial sensor attacks that verifying the $2s$-sparse observability condition (\cite{fawzi2014secure,shoukry2015event}) is NP-hard. Nevertheless,  polynomial time algorithms are provided to solve this problem, respectively, on a numerical system with bounded maximum geometric multiplicities, and on a structured system with bounded matching deficiencies, which are often met by practical engineering systems.  Numerical experiments show that when we have no access to the exact values of system matrices, the structured system based algorithms could be alternative for observability robustness evaluation. %An interesting extension  is to optimize some frequently-used quantitative observability metrics under cardinality-constrained sensor removals, such as trace of the inverse observability Gramian. %
Two problems, namely, the minSRO and the minSBO, are considered in this paper to study the observability robustness under sensor failures. It is shown that all these problems are NP-hard,  even restricted to some special cases, both for a numerical system and a structured one, thus confirming the conjecture that verifying the $2s$-sparse observability condition in \cite{fawzi2014secure,shoukry2015event} is NP-hard. The minSROs, both for a numerical system and a structured one, have a CCSM structure that is hard to approximate in general.  Nevertheless, under the constant bound assumption on the system eigenvalue geometric multiplicities or the matching deficiencies, by leveraging the rank-one update property of the rank function, polynomial time ({in the system dimensions and number of sensors}) exact algorithms are given for the minSROs. Numerical experiments on the minSRO show that
when we have no access to the exact values of system matrices,
the structured system based algorithms could be alternative for the
observability robustness evaluation.  Currently, we are not able to give (approximation) algorithms for the general minSBOs, which is left for future research.   %For space limitations,  some corollaries of results in this paper, an acceleration technique of Algorithm 2 based on  Dulmage-Mendelsohn decomposition, as well as some numerical simulations are provided in the extended version of this paper (see the updated version of \cite{zhang2018arxiv}).

%The problem of determining the minimal number of sensors whose removal destroys system observability is considered. This problem is closely related to secure state estimation under adversarial sensor attacks, and its dual is also the inverse to the MCP \cite{A.Ol2014Minimal, T2016On}. It is shown that this problem is NP-hard, both for a numerically specific system, and for a structured system, thus confirming the conjecture that verifying the $2s$-sparse observability condition in \cite{fawzi2014secure,shoukry2015event} is NP-hard. Nevertheless, polynomial time algorithms are provided to solve this problem, respectively, on a numerical system with bounded maximum geometric multiplicities, and on a structured system with bounded matching deficiencies, which are often met by practical engineering systems.  Numerical experiments show that when we have no access to the exact values of system matrices, the structured system based algorithms could be alternative for observability robustness evaluation. For space limitations, we have provided some corollaries of results in this paper, an acceleration technique of Algorithm 2 based on  Dulmage-Mendelsohn decomposition, as well as some numerical simulations in the extended version of this paper (see the updated version of \cite{zhang2018arxiv}).

%Although most techniques/algorithms in this paper are designed for observabilitys, they may be applied to some other fields, like the minimal cost critical columns in the linear degeneracy problem \cite{Khachiyan_1995}.

%\bibliographystyle{unsrt}
%\bibliography{References}

\bibliographystyle{elsarticle-num}
{\footnotesize
%\bibliography{References}
\bibliography{yuanz3}
}

\end{document}